\documentclass[12pt]{amsproc}

\usepackage{amsmath,amssymb}

\copyrightinfo{2000}{American Mathematical Society}

\newcommand{\cA}{{\mathcal A}}
\newcommand{\cD}{{\mathcal D}}
\newcommand{\cH}{{\mathcal H}}

\newcommand{\bN}{{\mathbb N}}
\newcommand{\bQ}{{\mathbb Q}}
\newcommand{\bR}{{\mathbb R}}
\newcommand{\bC}{{\mathbb C}}

\newtheorem{definition}{definition}[section]
\newtheorem{lemma}{Lemma}[section]
\newtheorem{theorem}{Theorem}[section]
\newtheorem{proposition}{Proposition}[section]

\newtheorem{remark}{Remark}[section]

\numberwithin{equation}{section}

\begin{document}

\title[Associated and quasi associated homogeneous distributions]
{Associated and quasi associated homogeneous distributions
(generalized functions)}

\author{V.~M.~Shelkovich}
\address{Department of Mathematics, St.-Petersburg State Architecture
and Civil Engineering University, 2 Krasnoarmeiskaya 4, 190005,
St. Petersburg, Russia.}
\email{shelkv@vs1567.spb.edu}

\thanks{The author was supported in part by
DFG Project 436 RUS 113/809/0-1 and Grant
05-01-04002-NNIOa of Russian Foundation for Basic Research.}


\date{}


\begin{abstract}
In this paper analysis of the concept of {\it associated homogeneous
distributions} (generalized functions) is given, and some
problems related to these distributions are solved.
It is proved that (in the one-dimensional case) there exist {\it only\/}
{\it associated homogeneous distributions}  of order $k=1$.
Next, we introduce a definition of {\it quasi associated homogeneous
distributions} and provide a mathematical description of all quasi
associated homogeneous distributions and their Fourier transform.
It is proved that the class of {\it quasi associated homogeneous
distributions} coincides with the class of distributions introduced
by Gel$'$fand and Shilov~\cite[Ch.I,\S 4.]{G-Sh} as the class of
{\it associated homogeneous distributions}.
For the multidimensional case it is proved that $f$ is a
{\it quasi associated homogeneous distribution} if and only if
it satisfies the Euler type system of differential equations.
A new type of $\Gamma$-functions generated by quasi associated
homogeneous distributions is defined.
\end{abstract}

\maketitle

\setcounter{equation}{0}

\section{Introduction}
\label{s1}

\subsection{Associated homogeneous distributions.\/}\label{s1.1}
First, the concept of {\it associated homogeneous distribution
{\rm(}AHD{\rm)}} (for the one-dimensional case) was introduced
by I.~M.~Gel$'$fand and G.~E.~Shilov in the book~\cite[Ch.I,\S 4.1.]{G-Sh}.
Let us repeat their reasoning by almost exact quoting.

Let us define the dilatation operator on the space ${\cD}'(\bR)$
by the formula $U_{a}f(x)=f(ax)$, $a >0$.
The definition of a {\it homogeneous distribution {\rm(}HD{\rm)}}
is the following.
\begin{definition}
\label{de2} \rm
(~\cite[Ch.I,\S 3.11.,(1)]{G-Sh},~\cite[Ch.X,8.]{Vil},~\cite[3.2.]{Hor})
A distribution $f_0\in {\cD}'(\bR)$ is said to be {\it homogeneous}
of degree~$\lambda$ if for any $a >0$ and $\varphi \in {\cD}(\bR)$
we have
$$
\Bigl\langle f_0(x),\varphi\Big(\frac{x}{a}\Big) \Bigr\rangle
=a^{\lambda+1} \bigl\langle f_0(x),\varphi(x) \bigr\rangle,
$$
i.e.,
\begin{equation}
\label{1*}
U_{a}f_0(x)=f_0(ax)=a^{\lambda}f_0(x).
\end{equation}
\end{definition}

Thus a HD of degree~$\lambda$ is an eigenfunction of {\it any} dilatation
operator $U_{a}$, $a>0$ with the eigenvalue $a^{\lambda}$, where
$\lambda\in \bC$, and $\bC$ is the set of complex numbers.

It is well known that ``in addition to an eigenfunction belonging
to a given eigenvalue, a linear transformation will ordinarily
also have so-called {\it associated functions} of various
orders''~\cite[Ch.I,\S 4.1.]{G-Sh}. The functions $f_1,f_2,\dots,f_k,\dots$
are said to be {\it associated with the eigenfunction} $f_0$ of the
transformation $U$ if
\begin{equation}
\label{1}
\begin{array}{rcl}
\displaystyle
Uf_0&=&cf_0, \\
\displaystyle
Uf_1&=&cf_1+df_0, \\
\displaystyle
Uf_2&=&cf_2+df_1, \\
\displaystyle
\cdots&\cdot&\cdots\cdots\cdots\cdot, \\
\displaystyle
Uf_k&=&cf_k+df_{k-1}, \\
\displaystyle
\cdots&\cdot&\cdots\cdots\cdots\cdot \\
\end{array}
\end{equation}
Consequently, $U$ reproduces an {\it associated function} of $k$th
order except for some multiple {\it associated function} of $(k-1)$th
order.

{\bf (D1)} Taking these facts into account, in the book~\cite[Ch.I,\S 4.1.]{G-Sh},
by analogy with Definition (\ref{1}) the following definition is
introduced: a function $f_1(x)$ is said to be {\it associated
homogeneous} of order $1$ and of degree~$\lambda$ if for any $a >0$
\begin{equation}
\label{1.1}
f_1(ax)=a^{\lambda}f_1(x)+h(a)f_{0}(x),
\end{equation}
where $f_{0}$ is a homogeneous function of degree $\lambda$. Here, in view of
(\ref{1*}) and (\ref{1}), $c=a^{\lambda}$. As the second step,
in~\cite[Ch.I,\S 4.1.]{G-Sh} it is proved that up to a constant factor
\begin{equation}
\label{1.2}
h(a)=a^{\lambda}\log{a}.
\end{equation}
Thus, by setting in the relation (\ref{1.1}) $c=a^{\lambda}$ and
$d=h(a)=a^{\lambda}\log{a}$, Definition~(\ref{1.1}) reads as the
following.

\begin{definition}
\label{de3} \rm
(Gel$'$fand and Shilov~\cite[Ch.I,\S4.1.,(1),(2)]{G-Sh})
A distribution $f_1\in {\cD}'(\bR)$ is called {\it associated
homogeneous distribution {\rm(}AHD{\rm)}} of order $1$ and of
degree~$\lambda$ if for any $a >0$ and $\varphi \in {\cD}(\bR)$
$$
\Bigl\langle f_1,\varphi\Big(\frac{x}{a}\Big) \Bigr\rangle
=a^{\lambda+1} \bigl\langle f_1,\varphi \bigr\rangle
+ a^{\lambda+1}\log{a}\bigl\langle f_{0},\varphi \bigr\rangle,
$$
i.e.,
$$
U_{a}f_1(x)=f_1(ax)=a^{\lambda}f_1(x)+a^{\lambda}\log{a}f_{0}(x),
$$
where $f_{0}$ is a homogeneous distribution of degree $\lambda$.
\end{definition}

It is clear that the class of AHDs of order $k=0$ coincides with
the class of HDs.

In the end, according to (\ref{1*}), (\ref{1}), (\ref{1.2}), using
$c=a^{\lambda}$ and $d=h(a)=a^{\lambda}\log{a}$, the following definition
is introduced.

\begin{definition}
\label{de4} \rm
(Gel$'$fand and Shilov~\cite[Ch.I,\S 4.1.,(3)]{G-Sh})
A distribution $f_k\in {\cD}'(\bR)$ is called an {\it AHD\/} of
order $k$, $k=2,3,\dots$ and of degree~$\lambda$ if for any $a >0$
and $\varphi \in {\cD}(\bR)$
\begin{equation}
\label{4}
\Bigl\langle f_k,\varphi\Big(\frac{x}{a}\Big) \Bigr\rangle
=a^{\lambda+1} \bigl\langle f_k,\varphi \bigr\rangle
+ a^{\lambda+1}\log{a}\bigl\langle f_{k-1},\varphi \bigr\rangle,
\end{equation}
where $f_{k-1}$ is an {\it AHD\/} of order $k-1$ and of degree~$\lambda$.
\end{definition}

In the book~\cite[Ch.I,\S 4]{G-Sh} (see also the paper~\cite[Ch.X,8.]{Vil})
it is stated (without proof) the following.

\begin{proposition}
\label{pr1}
Any AHD of order $k$ and of degree~$\lambda$ is a linear combination of
the following linearly independent AHDs of order $k$, \ $k=1,2,\dots$
and of degree~$\lambda$:

{}{\rm(a)} $x_{\pm}^{\lambda}\log^k x_{\pm}$ for $\lambda \ne -1,-2,\dots$;

{}{\rm(b)} $P\big(x_{\pm}^{-n}\log^{k-1}x_{\pm}\big)$ for $\lambda=-1,-2,\dots$;

{}{\rm(c)} $(x\pm i0)^{\lambda}\log^k(x\pm i0)$ for all $\lambda$.
\end{proposition}

Definitions of the above distributions are given in Sec.~\ref{s3}.

\subsection{Main results and contents of the paper.}\label{s1.2}
In this paper analysis of the concept of {\it associated homogeneous
distributions {\rm(}generalized functions{\it)}} is given, and some
problems related to this class of distributions are solved.

Unfortunately, as it follows from Sec.~\ref{s2}, Definition~\ref{de4}
(Gel$'$fand and Shilov) of {\it AHD\/} for $k\ge 2$ is
{\it self-contradictory}. In particular, it comes into conflict with
Proposition~\ref{pr1}.
In Sec.~\ref{s2}, we prove that an AHD of order $k$ is reproduced
by the dilatation operator $U_{a}$ (for all $a >0$) up to an AHD
of order $k-1$ {\it only} if $k=1$.
Thus in Definition~\ref{de4} the recursive step for $k=2$ is impossible.
Consequently, there exist {\it only\/} AHDs of order $k=0$, i.e., HDs
(given by Definition~\ref{de2}) and of order $k=1$ (given by
Definition~(\ref{1.1}) or Definition~\ref{de3}). Definition~\ref{de4}
(from the book~\cite[Ch.I,\S 4.1.,(3)]{G-Sh}), which defines AHDs
of order $k\ge 2$ describes an empty class.

The cause is the following: any HD is an eigenfunction
of {\it all} dilatation operators $U_{a}f_0(x)=f_0(ax)$ (for {\it all} $a>0$),
while any AHD is an eigenfunction of {\it all} dilatation operators
only for $k=1$.

In Sec.~\ref{s3}, we study the symmetry of the class of distributions
mentioned in Proposition~\ref{pr1} under the action of dilatation
operators $U_a$, $a>0$.

In Sec.~\ref{s4}, results of Sec.~\ref{s3} lead to a natural
generalization of the notion of the {\it associated eigenvector}
(\ref{1}) and imply our Definition~\ref{de6} of a {\it quasi associated
homogeneous distribution {\rm(}QAHD{\rm)}} of degree $\lambda$ and
of order $k$ by relation
$$
U_{a}f_k(x)=f_k(ax)=a^{\lambda}f_k(x)+\sum_{r=1}^{k}h_r(a) f_{k-r}(x),
\quad k=0,1,2,\dots, \quad  \forall \ a>0,
$$
where $f_{k-r}(x)$ is a QAHD of order $k-r$, \ $h_r(a)$ is a
differentiable function, $r=1,2,\dots,k$. (Here for $k=0$ we suppose
that the sum in the right-hand side of the last relation is empty.)

Thus the QAHD of order $k$ is reproduced by the dilatation operator
$U_{a}$ (for all $a >0$) up to a linear combination of QAHDs of orders
$k-1,k-2,\dots,0$ (see (\ref{50})).
Here the dilatation operator $U_{a}$ acts as a discrete convolution
$$
U_af_k(x)=f_k(ax)=\big(f(x)*h(a)\big)_{k}
$$
of sequences $f(x)=\{f_{0}(x),f_{1}(x),f_{2},\dots\}$ and
$h(a)=\{h_{0}(a)=a^{\lambda},h_{1}(a),h_{2}(a),\dots\}$.

$\bullet$ According to Theorem~\ref{th3}, in order to introduce a QAHD
of degree $\lambda$ and of order $k$ one can use
Definition~\ref{de7} instead of Definition~\ref{de6}, i.e., the relation
\begin{equation}
\label{15.0}
U_{a}f_k(x)=f_{k}(ax)=a^{\lambda}f_{k}(x)
+\sum_{r=1}^{k}a^{\lambda}\log^r{a} f_{k-r}(x),
\quad  \forall \ a>0.
\end{equation}
$k=0,1,2,\dots,$. Here for $k=0$ we suppose that the sum in the
right-hand side of the last relation is empty.

$\bullet$ By differentiating relation (\ref{15.0}), it easy to prove
by induction that if $f_{k}$ is a QAHD of degree $\lambda$ and of order
$k$ then its derivative $\frac{df_{k}}{d\lambda}$ with respect to
$\lambda$ is a QAHD of degree $\lambda$ and of order $k+1$.

$\bullet$ The sum of a QAHD of degree~$\lambda$ and of order~$k$,
and a QAHD of degree~$\lambda$ and of order~$r\le k-1$ is a QAHD
of degree~$\lambda$ and of order~$k$.

$\bullet$ In view of Definitions~\ref{de2},~\ref{de3},~\ref{de7}, the
classes of QAHDs of orders $k=0$ and $k=1$ coincide with the class
of HDs and the class of AHDs (in the Gel$'$fand and Shilov sense) of order
$k=1$, respectively.

$\bullet$ According to Theorems~\ref{th2},~\ref{th3}, the class of
all QAHDs coincides with the class of distributions
$$
{\cA\cH}_0(\bR)={\rm span}\{x_{\pm}^{\lambda}\log^k x_{\pm}, \
P\big(x_{\pm}^{-n}\log^{m-1}x_{\pm}\big):
\qquad\qquad\qquad\qquad\qquad\qquad
$$
$$
\qquad\qquad\qquad\qquad
\lambda \ne -1,-2,\dots,-n,\dots; \quad n, m \in \bN, \, k\in \{0\}\cup\bN \},
$$
introduced in the Gel$'$fand and Shilov book~\cite[Ch.I,\S 4.]{G-Sh}
as the class of AHDs (see Proposition~\ref{pr1}).

$\bullet$ According to Lemma~\ref{lem2}, QAHDs of different degrees and
orders are linear independent.

In Sec.~\ref{s5}, multidimensional QAHDs are
introduced. By Theorem~\ref{th6} it is proved that $f_k(x)$ is a QAHD
of order $k$, $k\ge 1$ if and only if it satisfies
the Euler type system of differential equations.
This result generalizes the well-known classical statement for homogeneous
distributions (see Theorem~\ref{th5}).

In Sec.~\ref{s6}, a mathematical description of the Fourier transform
of QAHDs is given for the multidimensional case.
Moreover, $\Gamma$-functions of a new type generated by QAHDs are defined.
In particular, for $k=1$ these $\Gamma$-functions are calculated and thweir
properties derived.

\begin{remark}
\label{rem1.1} \rm
In the papers~\cite{Al-Kh-Sh1},~\cite{Al-Kh-Sh2} a definition of an
{\it associated homogeneous $p$-adic distribution} was introduced and
mathematical description of all associated homogeneous distributions
and their Fourier transform was provided. This definition is the following:
$f\in {\cD}'(\bQ_p)$ is an {\it associated homogeneous distribution\/}
of degree~$\pi_{\alpha}(x)$ and order~$k$, \ $k=1,2,3\dots$, if for all
$\varphi \in {\cD}(\bQ_p)$ and $t \in \bQ_p^*$
\begin{equation}
\label{15}
\Bigl\langle f,\varphi\Big(\frac{x}{t}\Big) \Bigr\rangle
=\pi_{\alpha}(t)|t|_p \bigl\langle f,\varphi \bigr\rangle
+\sum_{j=1}^{k}\pi_{\alpha}(t)|t|_p\log_p^j|t|_p
\bigl\langle f_{k-j},\varphi \bigr\rangle
\end{equation}
where $f_{k-j}$ is an associated homogeneous distribution of
degree~$\pi_{\alpha}(x)$ and order $k-j$, \ $j=1,2,\dots,k$, i.e.
\begin{equation}
\label{15.1}
f(tx)=\pi_{\alpha}(t)f(x)
+\sum_{j=1}^{k}\pi_{\alpha}(t)\log_p^j|t|_pf_{k-j}(x),
\quad t \in \bQ_p^*.
\end{equation}
Here $\bQ_p$ is the field of $p$-adic numbers, $\bQ_p^*=\bQ_p\setminus \{0\}$
is its multiplicative group; $\pi_{\alpha}$ is a multiplicative character
of the field $\bQ_p$; ${\cD}(\bQ_p)$ is the linear space of
locally-constant $\bC$-value functions on $\bQ_p$ with compact supports,
${\cD}'(\bQ_p)$ is the set of all linear functionals on ${\cD}(\bQ_p)$.

One can see that a ``correct'' Definition~\ref{de7} of a
{\it quasi associated homogeneous distribution\/} is adaptation
(to the case $\bR$ instead of the field $\bQ_p$) of Definition~(\ref{15}),
(\ref{15.1}). However, in~\cite{Al-Kh-Sh1},~\cite{Al-Kh-Sh2} $p$-adic analog
of Theorem~\ref{th3} has not been proved.
\end{remark}

\section{Historical background, analysis, and comments.\/}
\label{s2}

{\bf (D2)} In contradiction to Definition~\ref{de4} (Gel$'$fand and
Shilov), in the paper of N.~Ya.~Vilenkin~\cite[Ch.X,8.]{Vil},
based on the book~\cite{G-Sh}, the following definition is used.
\begin{definition}
\label{de4-V} \rm
(Vilenkin~\cite[Ch.X,8.]{Vil})
A distribution $f_k\in {\cD}'(\bR)$ is called an {\it AHD\/} of order
$k$, $k=2,3,\dots$ and of degree~$\lambda$ if for any $a >0$
$$
f_k(ax)=a^{\lambda}f_k(x)+a^{\lambda}\log^k{a}f_{k-1}(x),
$$
where $f_{k-1}$ is an {\it AHD\/} of order $k-1$ and of degree~$\lambda$.
\end{definition}

Here an analog of relation (\ref{4}) is used, where in the right-hand
side of (\ref{4}) the term $\log{a}$ is replaced by $\log^k{a}$.

In the paper~\cite[Ch.X,8.]{Vil} Proposition~\ref{pr1} is also given
(without proof).

{\bf Comments.} (i) For example, according to Proposition~\ref{pr1},
$\log^2(x_{\pm})$ is an AHD of order $2$ and of degree~$0$. Nevertheless,
we have for all $a>0$
$$
\log^2(ax_{\pm})=\log^2 x_{\pm}+2\log a \log x_{\pm}+\log^2a.
$$
which contradicts Definitions~\ref{de4},~\ref{de4-V}.

(ii) In Sec.~\ref{s3}, relations (\ref{11}), (\ref{16}) imply that
in compliance with Proposition~\ref{pr1}, $x_{\pm}^{\lambda}\log x_{\pm}$
and $P\big(x_{\pm}^{-n}\big)$ are AHDs of order $k=1$ and of degree
$\lambda$ and $-n$, respectively (in the sense of Definition~\ref{de3}).
However, for $k\ge 2$, relations (\ref{11}), (\ref{16}) imply that
$x_{\pm}^{\lambda}\log^k x_{\pm}$ and $P\big(x_{\pm}^{-n}\log^{k-1} x_{\pm}\big)$
{\it are not AHDs} of order $k$ (in the sense of the above
Definition~\ref{de4} or Definition\ref{de4-V}). This contradicts
to Proposition~\ref{pr1}.

(iii) It remains to note that the assumption that an AHD of degree~$\lambda$
and of order $k$, \ $k\ge 2$ is defined by the Gel$'$fand--Shilov
Definition~\ref{de4}, contradicts some results on {\it distributional
quasi-asymptotics\/}. Indeed, if we temporarily assume that an AHD
of degree~$\lambda$ and of order $k$ is defined by Definition~\ref{de4},
in view of (\ref{4}), we have the asymptotic formulas:
$$
\begin{array}{rcl}
\displaystyle
f_{k}(ax)&=&a^{\lambda}f_{k}(x)+a^{\lambda}\log{a} f_{k-1}(x),
\qquad a \to \infty, \medskip \\
\displaystyle
f_{k}\Big(\frac{x}{a}\Big)&=&a^{-\lambda}f_{k}(x)-a^{-\lambda}\log{a} f_{k-1}(x),
\qquad a \to \infty.
\end{array}
$$
Here the coefficients of the {\it leading term\/} of both asymptotics
$f_{k-1}(x)$ and $-f_{k-1}(x)$ are AHDs of degree $\lambda$ and of order $k-1$.

In view of the above asymptotics, and according to~\cite{D-Zav1},
~\cite[Ch.I,Sec.~3.3.,Sec.~3.4.]{Vl-D-Zav}, the distribution $f_{k}$
has the {\it distributional quasi-asymptotics\/} $f_{k-1}(x)$ at
infinity with respect to an {\it automodel\/} function $a^{\lambda}\log{a}$,
and the distributional {\it quasi-asymptotics\/} $-f_{k-1}(x)$ at zero
with respect to an {\it automodel\/} function $a^{-\lambda}\log^k{a}$:
\begin{equation}
\label{9*}
\begin{array}{rclllll}
\displaystyle
f_{k}(x) &\stackrel{{\cD}'}{\sim}& f_{k-1}(x), &&x \to \infty
&& \big(a^{\lambda}\log{a}\big), \smallskip \\
\displaystyle
f_{k}(x) &\stackrel{{\cD}'}{\sim}& -f_{k-1}(x), &&x \to 0
&& \big(a^{-\lambda}\log{a}\big).
\end{array}
\end{equation}
Here both distributional quasi-asymptotics are AHDs of degree $\lambda$ and
of order $k-1$ (in the sense of Definition~\ref{de4}), \ $k\ge 2$. However,
according to~\cite{D-Zav1}, ~\cite[Ch.I,Sec.~3.3.,Sec.~3.4.]{Vl-D-Zav},
a distributional quasi-asymptotics is a {\it homogeneous distribution}.
Thus we have a contradiction.

\begin{remark}
\label{rem3} \rm
Let $f_{k}\in {\cA\cH}_0(\bR)$ be a QAHD of
degree~$\lambda$ and of order $k$, \ $k\ge 1$. In view of Definition~\ref{de7}
(see (\ref{15.0})), we have the asymptotic formulas:
\begin{equation}
\label{69}
\begin{array}{rcl}
\displaystyle
f_{k}(ax)&=&a^{\lambda}f_{k}(x)+\sum_{r=1}^{k}a^{\lambda}\log^r{a} f_{k-r}(x),
\qquad a \to \infty, \medskip \\
\displaystyle
f_{k}\Big(\frac{x}{a}\Big)&=&a^{-\lambda}f_{k}(x)
+\sum_{r=1}^{k}(-1)^{r}a^{-\lambda}\log^r{a} f_{k-r}(x),
\qquad a \to \infty.
\end{array}
\end{equation}
Here the coefficients of the {\it leading term\/} of both
asymptotics are homogeneous distributions $f_0$ and $(-1)^k f_0$
of degree~$\lambda$.

According to~\cite{D-Zav1},~\cite[Ch.I,Sec.~3.3.,Sec.~3.4.]{Vl-D-Zav}
and formulas (\ref{69}), the distribution $f_{k}$ has the {\it distributional
quasi-asymptotics\/} $f_0(x)$ at infinity with respect to an {\it automodel\/}
function $a^{\lambda}\log^k{a}$, and the distributional {\it quasi-asymptotics\/}
$(-1)^k f_0(x)$ at zero with respect to an {\it automodel\/} function
$a^{-\lambda}\log^k{a}$:
\begin{equation}
\label{70}
\begin{array}{rclllll}
\displaystyle
f_{k}(x) &\stackrel{{\cD}'}{\sim}& f_0(x), &&x \to \infty
&& \big(a^{\lambda}\log^k{a}\big), \smallskip \\
\displaystyle
f_{k}(x) &\stackrel{{\cD}'}{\sim}& (-1)^k f_0(x), &&x \to 0
&& \big(a^{-\lambda}\log^k{a}\big).
\end{array}
\end{equation}

In contrast to (\ref{9*}), both distributional quasi-asymptotics
(\ref{70}) are {\it homogeneous distributions\/}. This {\it is in compliance}
with the corresponding result from~\cite{D-Zav1},
~\cite[Ch.I,Sec.~3.3.,Sec.~3.4.]{Vl-D-Zav}: a distributional
quasi-asymptotics is a {\it homogeneous distribution}.
Thus our Definition~\ref{de7}, unlike Definition~\ref{de4},
implies the ``correct'' results on distributional quasi-asymptotics.
\end{remark}

(iv) Let us make an attempt ``to preserve'' Definition~(\ref{de4})
by some minor technical modifications.

By analogy with relation (\ref{1.1}) we will seek a function
$h_1(a)$ such that if $f_2(x)$ is an AHD of order $2$ and of
degree~$\lambda$ then for any $a >0$
\begin{equation}
\label{9}
U_{a}f_2(x)=f_2(ax)=a^{\lambda}f_2(x)+h_1(a)f_{1}(x),
\end{equation}
where $f_{1}(x)$ is an AHD of order $1$ and of degree~$\lambda$.

Similarly to~\cite[Ch.I,\S 4.1.]{G-Sh}, using (\ref{9}) and Definition~\ref{de3},
we obtain
$$
f_2(abx)=(ab)^{\lambda}f_2(x)+h_1(ab)f_{1}(x)
=a^{\lambda}f_2(bx)+h_{1}(a)f_{1}(bx)
\qquad\qquad\quad
$$
$$
\quad
=a^{\lambda}\big(b^{\lambda}f_2(x)+h_{1}(b)f_{1}(x)\big)
+h_{1}(a)\big(b^{\lambda}f_1(x)+b^{\lambda}\log b {\widetilde f}_{0}(x)\big)
$$
$$
\qquad
=(ab)^{\lambda}f_2(x)
+\Big(a^{\lambda}h_{1}(b)+b^{\lambda}h_{1}(a)\Big)f_{1}(x)
+h_{1}(a)b^{\lambda}\log b {\widetilde f}_{0}(x),
$$
where ${\widetilde f}_{0}(x)$ is a HD of  degree~$\lambda$.
Then for all $a,b>0$:
$$
\big(h_1(ab)-a^{\lambda}h_{1}(b)+b^{\lambda}h_{1}(a)\big)f_{1}(x)
-h_{1}(a)b^{\lambda}\log b {\widetilde f}_{0}(x)=0.
$$

It is easy to prove that a HD of degree $\lambda$ and an AHD
of order $1$ and of degree~$\lambda$ are linear independent
(see below Lemma~\ref{lem2}).
Consequently, there are two possibilities. If $h_{1}(a)\equiv 0$
then, according to (\ref{9}), $f_2(x)$ is a HD of degree $\lambda$.
If ${\widetilde f}_{0}(x)\equiv 0$ then
$h_1(ab)=a^{\lambda}h_{1}(b)+b^{\lambda}h_{1}(a)$, $h_1(1)=0$.
As mentioned above, the last equation has solution (\ref{1.2}),
and, consequently, $f_2(x)$ is an AHD of order $1$ and of degree~$\lambda$.

Thus it is impossible even for $k=2$ to preserve relation (\ref{1})
for all dilatation operators $U_{a}f(x)=f(ax)$, $a >0$. Consequently,
{\it it is impossible\/} to construct an AHD of order $k\ge 2$
defined by relation (\ref{1}) with the coefficients $c=a^{\lambda}$ and
$d=h(a)=a^{\lambda}\log{a}$.

\begin{remark}
\label{rem1} \rm
Definitions~\ref{de3},~\ref{de4} are given in compliance with the
book~\cite[Ch.I,\S 4.1.,(3)]{G-Sh}. Thus, in the case of Definition~\ref{de3}
(which defines an AHD of order $1$) one can clearly see that a
distribution $f_{0}$ does not depend on $a$.
In the case of Definition~\ref{de4} (which defines an AHD of order $k$
for $k\ge 2$), there is no clearness about independence of $f_{k-1}$ from $a$.
However, it is impossible ``to preserve'' the definition~\cite[Ch.I,\S 4.1.]{G-Sh}
even if we suppose that a distribution $f_{k-1}$ may depend on the
variable $a$.

Indeed, if we suppose that in Definition~\ref{de4} $f_{k-1}$
may depend on $a$, we will need to define AHD of degree $\lambda$
and of order $k\ge 2$ by the following relation
$$
f_k(ax)=a^{\lambda}f_k(x)+e(a)f_{k-1}(x,a), \quad \forall \ a>0,
$$
where $f_{k-1}(x,a)$ is an AHD (with respect of $x$) of degree $\lambda$
and of order $k-1$. It is clear that it is {\it impossible to determine}
a function $e(a)$.
\end{remark}

Thus, Definition~\ref{de4} (from the book~\cite[Ch.I,\S 4.1.,(3)]{G-Sh}
as well as Definition~\ref{de4-V} (from the paper~\cite[Ch.X,8.]{Vil})
define an {\it empty class}, and, consequently, the recursive step for
$k=2$ is impossible.

{\bf (D3)} In the books of R.~Estrada and R.~P.~Kanwal~\cite{E-K1},
~\cite{E-K2}, according to (\ref{1*}), (\ref{1}), a concept of
an {\it associated homogeneous distribution\/} is defined recursively.

\begin{definition}
\label{de4-E-K} \rm
(Estrada and Kanwal~\cite[(2.6.19)]{E-K1},~\cite[(2.110)]{E-K2})
An {\it associated homogeneous distribution\/} $f_k\in {\cD}'(\bR)$ of
order $k$ and of degree~$\lambda$ is such that for any $a >0$
\begin{equation}
\label{4-E-K}
f_k(ax)=a^{\lambda}f_k(x)+ a^{\lambda}e(a)f_{k-1}(x),
\end{equation}
where $f_{k-1}$ is an {\it associated homogeneous distribution\/} of
order $k-1$ and of degree~$\lambda$, and $e(a)$ is some function.
\end{definition}

Next, in these books {\it it is stated} that formula (\ref{4-E-K})
(i.e., formula~\cite[(2.6.19)]{E-K1},~\cite[(2.110)]{E-K2})
{\it implies} the relation
\begin{equation}
\label{5-E-K}
e(ab)=e(a)+e(b),
\end{equation}
i.e.,
\begin{equation}
\label{6-E-K}
e(a)=K\log{a}
\end{equation}
for some constant $K$, which can be absorbed in $f_{k-1}$~\cite[p.67]{E-K1},
~\cite[p.76]{E-K2}.
Finally, the authors of these books conclude that in view of
(\ref{4-E-K})--(\ref{6-E-K}) one can define an {\it associated
homogeneous distribution\/} of order $k-1$ and of degree~$\lambda$ by
the following equality
\begin{equation}
\label{7-E-K}
f_k(ax)=a^{\lambda}f_k(x)+a^{\lambda}\log{a}f_{k-1}(x),
\quad \forall \, a>0,
\end{equation}
where $f_{k-1}$ is an {\it associated homogeneous distribution\/} of
order $k-1$ and of degree~$\lambda$.
Thus, Definition (\ref{7-E-K}) (Estrada and Kanwal) coincides
with Definition~\ref{de4} (Gel$'$fand and Shilov).

{\bf  Comments.} Let us prove that formula (\ref{4-E-K})
(i.e., formula~\cite[(2.6.19)]{E-K1},~\cite[(2.110)]{E-K2})
{\it does not imply} relation (\ref{5-E-K}) for any $k\ge 2$.
Indeed, in view of (\ref{4-E-K}) we have for any $a,b>0$
\begin{equation}
\label{8-E-K}
f_k(abx)=(ab)^{\lambda}f_k(x)+ (ab)^{\lambda}e(ab)f_{k-1}(x)
=a^{\lambda}f_k(bx)+ a^{\lambda}e(a)f_{k-1}(bx),
\end{equation}
and
\begin{equation}
\label{9-E-K}
\begin{array}{rcl}
\displaystyle
f_k(bx)=b^{\lambda}f_k(x)+ b^{\lambda}e(b)f_{k-1}(x), \smallskip \\
\displaystyle
f_{k-1}(bx)=b^{\lambda}f_{k-1}(x)+ b^{\lambda}e(b)f_{k-2}(x), \\
\end{array}
\end{equation}
where $f_{k-1}$ and $f_{k-2}$ are AHDs of degree~$\lambda$ and of
order $k-1$ and $k-2$, respectively, $e(a)$ is some function.
By substituting relations (\ref{9-E-K}) into~(\ref{8-E-K}), we obtain
$$
(ab)^{\lambda}f_k(x)+ (ab)^{\lambda}e(ab)f_{k-1}(x)
\qquad\qquad\qquad\qquad\qquad\qquad\qquad\qquad\qquad
$$
$$
=a^{\lambda}\big(b^{\lambda}f_k(x)+ b^{\lambda}e(b)f_{k-1}(x)\big)
+ a^{\lambda}e(a)\big(b^{\lambda}f_{k-1}(x)+ b^{\lambda}e(b)f_{k-2}(x)\big).
$$
Thus we have for all $a,b >0$
\begin{equation}
\label{10-E-K}
\big(e(ab)-e(a)-e(b)\big)f_{k-1}(x)-e(a)e(b)f_{k-2}(x)=0,
\end{equation}
$k=1,2,\dots$. Here we set $f_{-1}(x)=0$.

It is clearly seen that, in contrast to the above cited statement
from~\cite{E-K1},~\cite{E-K2}) relation (\ref{10-E-K}) is equivalent to relation
(\ref{5-E-K}) only if $f_{k-2}(x)=0$, i.e., $k=1$.

Indeed, setting $k=1$, we calculate that $e(a)=K\log{a}$, i.e.,
(\ref{7-E-K}) holds for $k=1$.

Let $k=2$. In this case using (\ref{10-E-K}) and (\ref{6-E-K}), we obtain
$$
\big(\log{ab}-\log{a}-\log{b}\big)f_{1}(x)-\log{a}\log{b}f_{0}(x)=0,
$$
i.e., $f_{0}(x)\equiv 0$, which means that $f_1(x)$ is
a {\it homogeneous} distribution, and consequently, we have a contradiction.

{\bf (D4)} It remains to note that in the book~\cite{Hor}, the concept of
AHD is not discussed. It is only stated that for the distribution
$P\big(x_{+}^{-n}\big)$ ``the homogeneity is partly lost''.
However, according to Definition~\ref{de3} and Proposition~\ref{pr1}
(Gel$'$fand and  Shilov) this distribution is AHD of
order $1$ and of degree~$-n$, i.e., has a special symmetry.

{\bf Conclusion.} The concept of {\it associated homogeneous function}
has a misty prehistory. According to the above result, a {\it direct transfer}
of the notion of the {\it associated eigenvector} to the case of distributions
{\it is impossible for $k\ge 2$}. This is connected to the fact that any HD
is an eigenfunction of {\it all} dilatation operators $U_{a}f(x)=f(ax)$
(for {\it all} $a>0$), while for $k\ge 2$ {\it no} distribution
$x_{\pm}^{\lambda}\log^k x_{\pm}$, $P\big(x_{\pm}^{-n}\log^{k-1}x_{\pm}\big)$
is an AHD of {\it all\/} the dilatation operators.

\section{Symmetry of the class of distributions ${\cA\cH}_0(\bR)$}
\label{s3}

The distributions mentioned in Proposition~\ref{pr1} (so-called
``pseudo-functions'') are defined as regularizations of slowly divergent
integrals.
So, for all $\varphi \in {\cD}(\bR)$ and for $Re \lambda>-1$ we set
\begin{equation}
\label{5}
\Bigl\langle x_{+}^{\lambda}\log^k x_{+},\varphi(x) \Bigr\rangle
\stackrel{def}{=}\int_{0}^{\infty}x^{\lambda}\log^k x \varphi(x)\,dx.
\end{equation}
For $Re \lambda>-n-1$, \ $\lambda\ne -1,-2,\dots,-n$, according
to~\cite[Ch.I,\S 4.2.,(2),(6)]{G-Sh}, we have
$$
\Bigl\langle x_{+}^{\lambda}\log^k x_{+},\varphi(x) \Bigr\rangle
=\int_{0}^{1}x^{\lambda}\log^k x \bigg(\varphi(x)
-\sum_{j=0}^{n-1}\frac{x^{j}}{j!}\varphi^{(j)}(0)\bigg)\,dx
\qquad\qquad
$$
\begin{equation}
\label{6}
\qquad\qquad
+\int_{1}^{\infty}x^{\lambda}\log^k x \varphi(x)\,dx
+\sum_{j=0}^{n-1}\frac{(-1)^k k!}{j!(\lambda+j+1)^{k+1}}\varphi^{(j)}(0).
\end{equation}
The last formula gives an analytical continuation of relation (\ref{5}).

The distribution $P\big(x_{+}^{-n}\log^k x_{+}\big)$ (is not a value of
distribution $x_{+}^{\lambda}\log^k x_{+}$ at the point $\lambda=-n$)
is the principal value of the function~$x_{+}^{-n}\log^k x_{+}$.
According to~\cite[Ch.I,\S 4.2.,(4),(7)]{G-Sh}, we have
$$
\Bigl\langle P\big(x_{+}^{-n}\log^k x_{+}\big),\varphi(x) \Bigr\rangle
\qquad\qquad\qquad\qquad\qquad\qquad\qquad\qquad\qquad\qquad\qquad
$$
\begin{equation}
\label{7}
\stackrel{def}{=}\int_{0}^{\infty}x^{-n}\log^k x \bigg(\varphi(x)
-\sum_{j=0}^{n-2}\frac{x^{j}}{j!}\varphi^{(j)}(0)
-\frac{x^{n-1}}{(n-1)!}\varphi^{(n-1)}(0)H(1-x)\bigg)\,dx
\end{equation}
where $H(x)$ is the Heaviside function.

Other distributions mentioned in Proposition~\ref{pr1} are defined
as the following.
\begin{equation}
\label{8}
\begin{array}{rcl}
\displaystyle
\Bigl\langle x_{-}^{\lambda}\log^k x_{-},\varphi(x)\Bigr\rangle
&\stackrel{def}{=}&
\Bigl\langle x_{+}^{\lambda}\log^k x_{+},\varphi(-x) \Bigr\rangle, \medskip \\
\displaystyle
\Bigl\langle P\big(x_{-}^{-n}\log^k x_{-}\big),\varphi(x)\Bigr\rangle
&\stackrel{def}{=}&
\Bigl\langle P\big(x_{+}^{-n}\log^k x_{+}\big),\varphi(-x)\Bigr\rangle. \\
\end{array}
\end{equation}
for all $\varphi \in {\cD}(\bR)$.

Distributions $(x\pm i0)^{\lambda}\log^k(x\pm i0)$ are represented as
linear combinations of distributions $x_{\pm}^{\lambda}\log^k x_{\pm}$,
$P\big(x_{\pm}^{-n}\log^k x_{\pm}\big)$~\cite[Ch.I,\S4.5.]{G-Sh}.
In particular, for all $\lambda$~\cite[Ch.I,\S3.6.]{G-Sh}
\begin{equation}
\label{8.1}
\begin{array}{rcl}
\displaystyle
(x\pm i0)^{\lambda}&=&x_{+}^{\lambda}+e^{\pm i\pi\lambda}x_{-}^{\lambda},
\quad \lambda \ne -n, \quad n\in \bN; \medskip \\
\displaystyle
(x\pm i0)^{-n}&=&P\big(x^{-n}\big)
\mp \frac{i\pi(-1)^{n-1}\delta^{n-1}(x)}{(n-1)!},
\end{array}
\end{equation}
where the distribution $P\big(x^{-n}\big)$ is called the principal
value of the function~$x^{-n}$. This distribution is a homogeneous
distribution of degree $-n$. The distribution
$(x\pm i0)^{\lambda}\log^k(x\pm i0)$ for $\lambda \ne -1,-2,\dots$
can be obtained by differentiating the first relation in (\ref{8.1})
with respect to $\lambda$.

Let us consider how distributions from the class ${\cA\cH}_0(\bR)$
(mentioned above in Proposition~\ref{pr1}) are transformed by
dilatation operators $U_a$, $a>0$.

{\bf 1.} For $Re \lambda>-1$, $k\in \bN$ and for all
$\varphi(x) \in {\cD}(\bR)$, \ $a>0$ definition (\ref{5}) implies
$$
\Bigl\langle x_{+}^{\lambda}\log^k x_{+},\varphi\Big(\frac{x}{a}\Big) \Bigr\rangle
=a^{\lambda+1}\int_{0}^{\infty}\xi^{\lambda}\log^k (a\xi)\varphi(\xi)\,d\xi
\qquad\qquad\qquad\qquad\qquad\qquad
$$
$$
\quad
=a^{\lambda+1}\sum_{j=0}^{k}\log^j a C_{k}^{j}
\int_{0}^{\infty}\xi^{\lambda}\log^{k-j}\xi \varphi(\xi)\,d\xi
$$
\begin{equation}
\label{10}
\qquad
=a^{\lambda+1}\bigl\langle x_{+}^{\lambda}\log^k x_{+},\varphi(x) \bigr\rangle
+\sum_{j=1}^{k}a^{\lambda+1}\log^j a
\bigl\langle f_{\lambda;k-j}(x),\varphi(x)\bigr\rangle,
\end{equation}
where $f_{\lambda;k-j}(x)=C_{k}^{j} x_{+}^{\lambda}\log^{k-j}x_{+}$,
$C_{k}^{j}$ are binomial coefficients, $j=1,2,\dots,k$.
For all $\lambda\ne -1,-2,\dots$ we define (\ref{10}) by means of analytic
continuation. Thus
\begin{equation}
\label{11}
\Bigl\langle x_{+}^{\lambda}\log^k x_{+}, \varphi\Big(\frac{x}{a}\Big)\Bigr\rangle
=a^{\lambda+1}\bigl\langle x_{+}^{\lambda}\log^k x_{+},\varphi(x) \bigr\rangle
\sum_{j=1}^{k}a^{\lambda+1}\log^j a
\bigl\langle f_{\lambda;k}(x),\varphi(x)\bigr\rangle,
\end{equation}
for all $\lambda\ne -1,-2,\dots$.

{\bf 2.} For $k\in \bN$ and for all $\varphi(x) \in {\cD}(\bR)$
definition (\ref{7}) implies the following relations.

(a) $0<a<1$:
$$
\Bigl\langle P\big(x_{+}^{-n}\log^k x_{+}\big),
\varphi\Big(\frac{x}{a}\Big)\Bigr\rangle
\qquad\qquad\qquad\qquad\qquad\qquad\qquad\qquad\qquad\qquad
$$
$$
=\int_{0}^{1}x^{-n}\log^k x \bigg(\varphi(x/a)
-\sum_{j=0}^{n-1}\frac{(x/a)^{j}}{j!}\varphi^{(j)}(0)\bigg)\,dx
\qquad\qquad
$$
$$
\qquad\qquad\qquad
+\int_{1}^{\infty}x^{-n}\log^k x \bigg(\varphi(x/a)
-\sum_{j=0}^{n-2}\frac{(x/a)^{j}}{j!}\varphi^{(j)}(0)\bigg)\,dx
$$
$$
\qquad
=a^{-n+1}\bigg\{\int_{0}^{1/a}\xi^{-n}\log^k(a\xi)\bigg(\varphi(\xi)
-\sum_{j=0}^{n-1}\frac{\xi^{j}}{j!}\varphi^{(j)}(0)\bigg)\,d\xi
\qquad\qquad
$$
$$
\qquad\qquad\qquad
+\int_{1/a}^{\infty}\xi^{-n}\log^k(a\xi) \bigg(\varphi(\xi)
-\sum_{j=0}^{n-2}\frac{\xi^{j}}{j!}\varphi^{(j)}(0)\bigg)\,d\xi \bigg\}
$$
$$
\qquad
=a^{-n+1}\bigg\{\int_{0}^{1}\xi^{-n}\log^k(a\xi)\bigg(\varphi(\xi)
-\sum_{j=0}^{n-1}\frac{\xi^{j}}{j!}\varphi^{(j)}(0)\bigg)\,d\xi
$$
$$
+\int_{1}^{\infty}\xi^{-n}\log^k(a\xi) \bigg(\varphi(\xi)
-\sum_{j=0}^{n-2}\frac{\xi^{j}}{j!}\varphi^{(j)}(0)\bigg)\,d\xi
-\frac{\varphi^{(n-1)}(0)}{(n-1)!}I_1\bigg\}
$$
\begin{equation}
\label{12}
=a^{-n+1}\bigg\{\sum_{r=0}^k\log^r a C_k^r
\langle P\big(x_{+}^{-n}\log^{k-r} x_{+}\big),\varphi(x) \rangle
-\frac{\varphi^{(n-1)}(0)}{(n-1)!}I_1\bigg\},
\end{equation}
where
$$
I_1=\int_{1}^{1/a}\frac{\log^k(a\xi)}{\xi}\,d\xi
=\sum_{r=0}^k\log^r a C_k^r\int_{1}^{1/a}\frac{\log^{k-r} \xi}{\xi}\,d\xi
\qquad\qquad\qquad\qquad\qquad
$$
\begin{equation}
\label{13}
\qquad\qquad\qquad
=\log^{k+1}a\sum_{r=0}^k C_k^r \frac{(-1)^{k+1-r}}{k+1-r}
=-\frac{1}{k+1}\log^{k+1}a.
\end{equation}

(b) $a=1$:
$$
\Bigl\langle P\big(x_{+}^{-n}\log^k x_{+}\big),
\varphi\Big(\frac{x}{a}\Big)\Bigr\rangle
=\bigl\langle P\big(x_{+}^{-n}\log^k x_{+}\big),\varphi(x) \bigr\rangle.
$$

(c) $a>1$:
$$
\Bigl\langle P\big(x_{+}^{-n}\log^k x_{+}\big),
\varphi\Big(\frac{x}{a}\Big)\Bigr\rangle
\qquad\qquad\qquad\qquad\qquad\qquad\qquad\qquad\qquad
$$
$$
=a^{-n+1}\bigg\{\int_{0}^{1/a}\xi^{-n}\log^k(a\xi)\bigg(\varphi(\xi)
-\sum_{j=0}^{n-1}\frac{\xi^{j}}{j!}\varphi^{(j)}(0)\bigg)\,d\xi
\qquad
$$
$$
\qquad\qquad\qquad
+\int_{1/a}^{\infty}\xi^{-n}\log^k(a\xi) \bigg(\varphi(\xi)
-\sum_{j=0}^{n-2}\frac{\xi^{j}}{j!}\varphi^{(j)}(0)\bigg)\,d\xi \bigg\}
$$
\begin{equation}
\label{14}
=a^{-n+1}\bigg\{\sum_{r=0}^k\log^r a C_k^r
\langle P\big(x_{+}^{-n}\log^{k-r} x_{+}\big),\varphi(x) \rangle
+\frac{\varphi^{(n-1)}(0)}{(n-1)!}I_2\bigg\},
\end{equation}
where
$$
I_2=\int_{1/a}^{1}\frac{\log^k(a\xi)}{\xi}\,d\xi=-I_1=\frac{1}{k+1}\log^{k+1}a.
$$

Thus, (\ref{12})--(\ref{14}) imply
$$
\Bigl\langle P\big(x_{+}^{-n}\log^k x_{+}\big),
\varphi\Big(\frac{x}{a}\Big)\Bigr\rangle
\qquad\qquad\qquad\qquad\qquad\qquad\qquad\qquad\qquad\qquad
$$
\begin{equation}
\label{16}
=a^{-n+1}\bigl\langle P\big(x_{+}^{-n}\log^k x_{+}\big),\varphi(x)\bigr\rangle
+\sum_{r=1}^{k+1} a^{-n+1}\log^r a
\bigl\langle f_{-n;k+1-r}(x),\varphi(x) \bigr\rangle
\end{equation}
where $f_{-n;0}(x)=\frac{(-1)^{n-1}}{(k+1)(n-1)!}\delta^{(n-1)}(x)$ \
and $f_{-n;k+1-r}(x)=C_k^r P\big(x_{+}^{-n}\log^{k-r} x_{+}\big)$, \
$r=1,2,\dots,k$.

For distributions $x_{-}^{\lambda}\log^k x_{-}$, \
$P\big(x_{-}^{-n}\log^k x_{-}\big)$ relations of the type (\ref{11}),
(\ref{16}) can be obtained from (\ref{8}).

\section{Quasi associated homogeneous distributions}
\label{s4}

\subsection{A class of distributions ${\cA\cH}_1(\bR)$.}\label{s4.1}
In Sec.~\ref{s1}, it is recognized that the dilatation
operator $U_a$ for all $a>0$ does not reproduce a distribution
of order $k$ from ${\cA\cH}_0(\bR)$ with accuracy up to
a distribution of order $(k-1)$ from ${\cA\cH}_0(\bR)$.
Moreover, in Sec.~\ref{s3}, it is recognized that the
dilatation operator $U_a$ acts in ${\cA\cH}_0(\bR)$ according to
formulas (\ref{11}), (\ref{16}).
Now by analogy with transformation laws (\ref{11}), (\ref{16})
we introduce the following definition.

\begin{definition}
\label{de5} \rm
A distribution $f_{\lambda;k}\in {\cD}'(\bR)$ is called a {\it distribution
of degree~$\lambda$ and of order $k$}, $k=0,1,2,\dots$, if for any $a >0$
and $\varphi \in {\cD}(\bR)$
\begin{equation}
\label{17}
\Bigl\langle f_{\lambda;k}(x),\varphi\Big(\frac{x}{a}\Big) \Bigr\rangle
=a^{\lambda+1} \bigl\langle f_{\lambda;k}(x),\varphi(x) \bigr\rangle
+\sum_{r=1}^{k}
a^{\lambda+1}\log^r{a}\bigl\langle f_{\lambda;k-r}(x),\varphi(x) \bigr\rangle,
\end{equation}
i.e.,
\begin{equation}
\label{17.1}
U_{a}f_{\lambda;k}(x)=f_{\lambda;k}(ax)=a^{\lambda}f_{\lambda;k}(x)
+\sum_{r=1}^{k}a^{\lambda}\log^r{a} f_{\lambda;k-r}(x),
\end{equation}
where $f_{\lambda;k-r}(x)$ is a {\it distribution of degree $\lambda$ and
of order $k-r$}, \ $r=1,2,\dots,k$. Here for $k=0$ we suppose that
sums in the right-hand side of (\ref{17}), (\ref{17.1}) are empty.
\end{definition}

Let us denote by ${\cA\cH}_1(\bR)$ a linear span of all {\it distributions
$f_{\lambda;k}(x)\in {\cD}'(\bR)$ of order $k$ and degree $\lambda$},
$\lambda\in \bC$, $k=0,1,2,\dots$, defined by Definition~\ref{de5}.

In view of Definitions~\ref{de2},~\ref{de3},~\ref{de5},
a HD of degree~$\lambda$ is a {\it distribution of order $k=0$ and degree~$\lambda$},
and an AHD of order $1$ and degree~$\lambda$ is a {\it distribution of order $k=1$
and degree~$\lambda$}.
According to (\ref{11}), (\ref{16}),
$x_{\pm}^{\lambda}\log^k x_{\pm}$, and $P\big(x_{\pm}^{-n}\log^{k-1}x_{\pm}\big)$
are {\it distributions of order $k$ and of degree $\lambda$, and $-n$, respectively}.
Thus ${\cA\cH}_0(\bR) \subset {\cA\cH}_1(\bR)$.

\begin{remark}
\label{rem2} \rm
A sum of a distribution of degree~$\lambda$ and of order~$k$
(from ${\cA\cH}_1(\bR)$) and a distribution of degree~$\lambda$
and of order~$r\le k-1$ (from ${\cA\cH}_1(\bR)$) is a distribution
of degree~$\lambda$ and of order~$k$ (from ${\cA\cH}_1(\bR)$).
\end{remark}

\begin{lemma}
\label{lem2}
Distributions from ${\cA\cH}_1(\bR)$ of different degrees and orders are
linear independent.
\end{lemma}

\begin{proof}
This lemma is proved in the same way as the analogous result
on linear independent homogeneous distributions from~\cite[\S3.11.,4.]{G-Sh}.

Suppose that
$$
c_1 f^1(x)+\cdots+c_m f^m(x)=0,
$$
where $f^s(x)\in {\cA\cH}_1(\bR)$ is a distribution of degree~$\lambda$ and of
order~$k_s$, such that all $\lambda_s$ or $k_s$, \ $s=1,2,\dots,m$ are different.
Then, by Definition~\ref{de5}, for all $a>0$ and $\varphi(x) \in {\cD}(\bR)$:
$$
c_1a^{\lambda_1}\bigg(\Bigl\langle f^1(x),\varphi(x)\Bigr\rangle
+\sum_{r=1}^{k_1}\log^r{a}
\Bigl\langle f^1_{k_1-r}(x),\varphi(x)\Bigr\rangle\bigg)
\qquad\qquad\qquad\qquad\qquad
$$
$$
\qquad
+ \cdots + c_ma^{\lambda_m}\bigg(\Bigl\langle f^m(x),\varphi(x)\Bigr\rangle
+\sum_{r=1}^{k_m}\log^r{a}
\Bigl\langle f^m_{k_m-r}(x),\varphi(x)\Bigr\rangle\bigg)=0,
$$
where $f^s_{k_s-r}(x)\in {\cA\cH}_1(\bR)$ is a distribution of degree~$\lambda$
and of order~$(k_s-r)$, \ $r=1,2,\dots,k_s$, \ $s=1,2,\dots,m$.

If all $\lambda_s$ are different, by choosing different values $a$,
it is easy to see that, $c_s\equiv 0$, \ $s=1,2,\dots,m$.

If, for example, $\lambda_1=\lambda_2$ and $k_1> k_2$, then for all
$a>0$ and $\varphi(x) \in {\cD}(\bR)$ we have
$$
c_1\Big(\Bigl\langle f^1(x),\varphi(x)\Bigr\rangle
+\sum_{r=1}^{k_1}\log^r{a}
\Bigl\langle f^1_{k_1-r}(x),\varphi(x)\Bigr\rangle\Big)
\qquad\qquad\qquad\qquad\qquad
$$
$$
\qquad\qquad
+ c_2\Big(\langle f^2(x),\varphi(x)\rangle
+\sum_{r=1}^{k_2}\log^r{a} \langle f^2_{k_2-r}(x),\varphi(x)\rangle\Big)=0.
$$
The last relation implies that $c_1f^1_{k_1-r}(x)=0$, \
$r=k_2+1,\dots,k_1$, and, consequently, $c_1\equiv 0$. Consequently, $c_2\equiv 0$.
\end{proof}

\begin{theorem}
\label{th2}
Every distribution $f\in {\cA\cH}_1(\bR)$ of degree $\lambda$
and order $k\in \bN$ {\rm(}up to a distribution of order $\le k-1${\rm)}
is a sum of linear independent distributions

{\rm (a)} \ $C x_{\pm}^{\lambda}\log^k x_{\pm}$, if $\lambda \ne -1,-2,\dots$;

{\rm (b)} \ $CP\big(x_{\pm}^{-n}\log^{k-1} x_{\pm}\big)$, if
$\lambda=-n$, $n\in \bN$, where $C$ is a constant.

Thus ${\cA\cH}_1(\bR)={\cA\cH}_0(\bR)$, i.e., the class ${\cA\cH}_1(\bR)$
coincides with the Gel$'$fand and Shilov class ${\cA\cH}_0(\bR)$ from
Proposition~{\rm\ref{pr1}}.
\end{theorem}

\begin{proof}
We prove this theorem by induction.
(a) Let us consider the case $\lambda \ne -1,-2,\dots$.
According to Definitions~\ref{de3},~\ref{de5}, a distribution
$f_1\in {\cA\cH}_1(\bR)$ of degree $\lambda$ and order~$k=1$ is an
AHD of degree $\lambda$ and order~$k=1$, and for all $a>0$
satisfies the equation
\begin{equation}
\label{20}
f_1(ax)=a^{\lambda}f_1(x)+a^{\lambda}\log{a} f_{0}(x),
\end{equation}
where $f_{0}(x)$ is a HD of degree $\lambda$.
In view of Theorem~\cite[Ch.I,\S 3.11.]{G-Sh},
$f_{0}(x)=A_1x_{+}^{\lambda}+A_2x_{-}^{\lambda}$, where $A_1$, $A_2$
are constants. If we differentiate (\ref{20}) with respect to $a$
and set $a=1$, we obtain the differential equation
\begin{equation}
\label{21}
xf_1'(x)=\lambda f_1(x)+A_1x_{+}^{\lambda}+A_2x_{-}^{\lambda}.
\end{equation}
For $\pm x>0$ the last equation can be integrated in the ordinary sense.

Thus, for $x>0$ equation (\ref{21}) coincides with
the equation $xf_1'(x)=\lambda f_1(x)+A_1x_{+}^{\lambda}$.
Integrating this equation, we obtain
$f_1(x)=A_1x_{+}^{\lambda}\log x_{+}+B_1x_{+}^{\lambda}$,
where $B_1$ is a constant. Similarly one can prove that
$f_1(x)=A_2x_{-}^{\lambda}\log x_{-}+B_2x_{-}^{\lambda}$ for $x<0$.
Thus the distribution $g(x)=f_1(x)-A_1x_{+}^{\lambda}
-A_2x_{-}^{\lambda}-B_1x_{+}^{\lambda}-B_2x_{-}^{\lambda}$
satisfies equation (\ref{21}) being concentrated at the point
$x=0$. Therefore, $g(x)=\sum_{m=0}^MC_m\delta^{(m)}(x)$, where $C_1,\dots,C_M$
are constants.
However, since $\delta^{(m)}(x)$ is a HD of degree $-m-1$, in view of
Lemma~\ref{lem2} $g(x)=0$. Thus
$$
f_1(x)=A_1x_{+}^{\lambda}\log x_{+}+A_2x_{-}^{\lambda}\log x_{-}
+B_1x_{+}^{\lambda}+B_2x_{-}^{\lambda}.
$$
Consequently, up to a distribution
$B_1x_{+}^{\lambda}+B_2x_{-}^{\lambda}\in {\cA\cH}_1(\bR)$ of order $0$, we
have $f_1(x)=A_1x_{+}^{\lambda}\log x_{+}+A_2x_{-}^{\lambda}\log x_{-}$.

Let us assume that a distribution $f_{k-1}(x)\in {\cA\cH}_1(\bR)$
of degree $\lambda$ and order $(k-1)$ is represented in the form
of a linear combination
\begin{equation}
\label{22.1}
f_{k-1}(x)=\sum_{j=0}^{k-1}\Big(A_{1j}x_{+}^{\lambda}\log^{j} x_{+}
+A_{2j}x_{-}^{\lambda}\log^{j} x_{-}\Big).
\end{equation}

A distribution $f_k\in {\cA\cH}_1(\bR)$ of degree~$\lambda$ and of order~$k\ge 2$
satisfies (\ref{17.1}) for all $a>0$. By differentiating this equation
with respect to $a$ and setting $a=1$, we obtain
\begin{equation}
\label{22}
xf_{k}'(x)=\lambda f_{k}(x)+f_{k-1}(x).
\end{equation}
Taking into account (\ref{22.1}) and integrating (\ref{22}) for $x\ne 0$,
we calculate
$$
f_{k}(x)=\sum_{j=0}^{k-1}\bigg(\frac{A_{1j}}{j+1}x_{+}^{\lambda}\log^{j+1} x_{+}
+\frac{A_{2j}}{j+1}x_{-}^{\lambda}\log^{j+1} x_{-}\bigg)
+B_{1}x_{+}^{\lambda}+B_{2}x_{-}^{\lambda},
$$
where $B_1$, $B_2$ are constant. By repeating the above reasoning
we obtain that
$$
f_{k}(x)=\frac{A_{1 \ k-1}}{k}x_{+}^{\lambda}\log^{k} x_{+}
+\frac{A_{2 \ k-1}}{k}x_{-}^{\lambda}\log^{k} x_{-}
$$
up to distributions of degree $\lambda$ and of order $\le k-1$.

Hence, by induction the case (a) is proved.

The case (b), when $\lambda=-n$, $n\in \bN$, can be proved similarly
to the case (a).
\end{proof}

\subsection{QAHDs.}\label{s4.2}
Taking into account relations (\ref{11}), (\ref{16}), and by analogy with
(\ref{1.1}) we introduce the following definition.

\begin{definition}
\label{de6} \rm
A distribution $f_k\in {\cD}'(\bR)$ is called a {\it quasi associated
homogeneous distribution\/} of
degree~$\lambda$ and of order $k$, $k=0,1,2,3,\dots$ if for any $a >0$
and $\varphi \in {\cD}(\bR)$
$$
\Bigl\langle f_k(x),\varphi\Big(\frac{x}{a}\Big) \Bigr\rangle
=a^{\lambda+1} \langle f_k(x),\varphi(x) \rangle
+\sum_{r=1}^{k}h_r(a)\langle f_{k-r}(x),\varphi(x) \rangle,
$$
i.e.,
\begin{equation}
\label{50}
U_{a}f_k(x)=f_k(ax)=a^{\lambda}f_k(x)+\sum_{r=1}^{k}h_r(a) f_{k-r}(x),
\end{equation}
where $f_{k-r}(x)$ is a QAHD of degree $\lambda$ and of order $k-r$, \
$h_r(a)$ is a differentiable function, $r=1,2,\dots,k$.
Here for $k=0$ we suppose that sums in the right-hand sides of the above
relations are empty.
\end{definition}

Let us denote by ${\cA\cH}(\bR)$ a linear span of all {\it QAHDs\/}
of order $k$ and degree $\lambda$, $\lambda\in \bC$, $k=0,1,2,\dots$,
defined by Definition~\ref{de6}. In view of Definition~\ref{de5},
${\cA\cH}_1(\bR) \subset {\cA\cH}(\bR)$.

\begin{theorem}
\label{th3}
Any QAHD $f_k(x)$ of degree $\lambda$
and of order $k$, \ $k=0,1,\dots$ {\rm(}see Definition~{\rm \ref{de6}})
is a distribution of degree $\lambda$ and of order $k$ {\rm(}from
${\cA\cH}_0(\bR)${\rm)} {\rm(}see Definition~{\rm \ref{de5}} and
Theorem~{\rm \ref{th2}}{\rm)}, i.e., $f_k(x)$ satisfies relation
{\rm (\ref{17.1})}.

Thus ${\cA\cH}(\bR)={\cA\cH}_0(\bR)$, i.e., the class ${\cA\cH}(\bR)$
coincides with the Gel$'$fand and Shilov class ${\cA\cH}_0(\bR)$ from
Proposition~{\rm\ref{pr1}}.
\end{theorem}

\begin{proof}
We prove this theorem by induction.

1. For $k=1$ this theorem is proved in the book~\cite[Ch.I,\S 4.1.]{G-Sh}
(see also Subsec.~\ref{s1.1}).

2. If $k=2$, according to Definition~\ref{de6}, for a QAHD
$f_2(x)$ of degree $\lambda$ and of order $k=2$ we have
\begin{equation}
\label{51}
f_2(ax)=a^{\lambda}f_2(x)+h_{1}(a)f_{1}(x)+h_{2}(a)f_{0}(x),
\quad \forall \ a>0,
\end{equation}
where $f_{1}(x)$ is an AHD of degree $\lambda$ and of order $k=1$, \
$f_{0}(x)$ is a HD of degree $\lambda$, and $h_{1}(a)$, $h_{2}(a)$
are the desired functions.

Taking into account that
$f_{1}(bx)=b^{\lambda}f_1(x)+b^{\lambda}\log b f^{(1)}_{0}(x)$,
where $f^{(1)}_{0}(x)$ is a HD of degree $\lambda$, in view of
(\ref{51}) and Definition~\ref{de3}, we obtain for all $a,b>0$:
$$
f_2(abx)=(ab)^{\lambda}f_2(x)+h_{1}(ab)f_{1}(x)+h_{2}(ab)f_{0}(x)
\qquad\qquad\qquad\qquad
$$
$$
=a^{\lambda}f_2(bx)+h_{1}(a)f_{1}(bx)+h_{2}(a)f_{0}(bx)
\qquad\qquad\quad
$$
$$
=a^{\lambda}\Big(b^{\lambda}f_2(x)+h_{1}(b)f_{1}(x)+h_{2}(b)f_{0}(x)\Big)
\qquad\qquad
$$
$$
\quad
+h_{1}(a)\Big(b^{\lambda}f_1(x)+b^{\lambda}\log b f^{(1)}_{0}(x)\Big)
+h_{2}(a)b^{\lambda}f_{0}(x)
$$
$$
=(ab)^{\lambda}f_2(x)
+\Big(a^{\lambda}h_{1}(b)+b^{\lambda}h_{1}(a)\Big)f_{1}(x)
\qquad\qquad
$$
$$
\quad\qquad\quad
+\Big(a^{\lambda}h_{2}(b)+h_{2}(a)b^{\lambda}\Big)f_{0}(x)
+h_{1}(a)b^{\lambda}\log b f^{(1)}_{0}(x).
$$
Obviously, this implies that for all $a,b>0$
$$
\Big(h_{1}(ab)-a^{\lambda}h_{1}(b)-b^{\lambda}h_{1}(a)\Big)f_{1}(x)
\qquad\qquad\qquad\qquad\qquad\qquad\qquad\qquad
$$
\begin{equation}
\label{52}
+\Big(h_{2}(ab)-a^{\lambda}h_{2}(b)-b^{\lambda}h_{2}(a)\Big)f_{0}(x)
-h_{1}(a)b^{\lambda}\log b f^{(1)}_{0}(x)=0.
\end{equation}

According to~\cite[Ch.I,\S 3.11.]{G-Sh}, there are two {\it linear independent\/}
HDs of degree~$\lambda$, such that every HD is their linear combination. Thus
there are two possibilities: either $f^{(1)}_{0}(x)$ and $f_{0}(x)$ are linear
independent HDs, or $f^{(1)}_{0}(x)=Cf_{0}(x)$, where $C$ is a constant.

Thus in the first case, since in view of Lemma~\ref{lem2} a HD and an AHD
of order $1$ are linear independent, relation (\ref{52}) implies
$h_{1}(a)=0$ and $h_2(ab)=a^{\lambda}h_2(b)+b^{\lambda}h_2(a)$.
The  solution of the last equation constructed in~\cite[Ch.I,\S 4.1.]{G-Sh}
(see also Subsec.~\ref{s1.1}), is given by (\ref{1.2}), i.e.,
$h_2(a)=a^{\lambda}\log{a}$. Thus, relation (\ref{51}), Definition~\ref{de3},
and Theorem~\ref{th2} imply that $f_2(x)$ is an AHD of order $1$.
Consequently, we obtain a {\it trivial solution}.

In the second case, in view of Lemma~\ref{lem2}, $f_{0}(x)$ and
$f_{1}(x)$ are linear independent, and, consequently, relation (\ref{52})
implies the system of functional equations:
\begin{equation}
\label{53}
\begin{array}{rcl}
\displaystyle
h_1(ab)&=&a^{\lambda}h_1(b)+b^{\lambda}h_1(a), \medskip\\
\displaystyle
h_2(ab)&=&a^{\lambda}h_2(b)+h_2(a)b^{\lambda}
+Ch_1(a)b^{\lambda}\log b, \quad \forall \ a,b>0, \\
\end{array}
\end{equation}
where $h_1(1)=0$, $h_2(1)=0$. According to~\cite[Ch.I,\S 4.1.]{G-Sh} (see also
(\ref{1.2})), $h_1(a)=a^{\lambda}\log{a}$. Then the second equation in (\ref{53})
implies that
\begin{equation}
\label{54}
h_2(ab)=h_2(a)b^{\lambda}+a^{\lambda}h_2(b)+C(ab)^{\lambda}\log{a}\log{b}
\end{equation}
and, consequently, the function
${\widetilde h}_2(a)=\frac{h_2(a)}{a^{\lambda}}$ satisfies the equation
\begin{equation}
\label{55}
{\widetilde h}_2(ab)={\widetilde h}_2(a)+{\widetilde h}_2(b)
+C\log{a}\log{b}, \quad \forall \ a,b>0.
\end{equation}
Making the change of variables $\psi_2(z)={\widetilde h}_2(e^{z})$, where
$\psi_2(0)=0$ and $a=e^{\xi}$, $b=e^{\eta}$, we can see that (\ref{55})
can be rewritten as
\begin{equation}
\label{56}
\psi_2(\xi+\eta)=\psi_2(\xi)+\psi_2(\eta)+C\xi\eta, \quad \forall \ \xi,\eta.
\end{equation}

We will seek a solution of equation (\ref{56}) in the class of differentiable
functions. Differentiating relation (\ref{56}) with respect to $\eta$,
we obtain for all $\xi,\eta$
$$
\psi_2'(\xi+\eta)=\psi_2'(\eta)+C\xi, \qquad
\psi_2(0)=0.
$$
Setting $\eta=0$ in the last equation, we have the differential equation
$$
\psi_2'(\xi)=\psi_2'(0)+C\xi, \quad \psi_2(0)=0,
$$
whose solution has the form
$$
\psi_2(\xi)=\psi_2'(0)\xi+\frac{C}{2}\xi^2.
$$
Since $a=e^{\xi}$, then ${\widetilde h}_2(a)=A_2\log a+\frac{C}{2}\log^2 a$
and
\begin{equation}
\label{57}
h_2(a)=A_2a^{\lambda}\log a+\frac{C}{2}a^{\lambda}\log^2 a,
\end{equation}
where $A_2={\widetilde h}_2'(1)=h_2'(1)$ is a constant.

By substituting functions $h_1(a)$, $h_2(a)$ given by (\ref{1.2}), (\ref{57})
into (\ref{51}), we obtain
$$
f_2(ax)=a^{\lambda}f_2(x)+a^{\lambda}\log a f_{1}(x)
+\Big(A_2a^{\lambda}\log a+\frac{C}{2}a^{\lambda}\log^2 a\Big)f_{0}(x).
$$
The last relation can be rewritten in the desired form (\ref{17.1}):
\begin{equation}
\label{58}
f_2(ax)=a^{\lambda}f_2(x)+a^{\lambda}\log a {\widetilde f}_{1}(x)
+a^{\lambda}\log^2 a{\widetilde f}_{0}(x), \quad \forall \ a>0,
\end{equation}
where ${\widetilde f}_{1}(x)=f_{1}(x)+A_2f_{0}(x)$ is
an AHD of degree $\lambda$ and of order $k=1$, \
${\widetilde f}_{0}(x)=\frac{C}{2}f_{0}(x)$ is a HD.
Thus $f_2(x)$ is a distribution of degree $\lambda$
and of order $2$ (in the sense of Definition~\ref{de5}),
and, according to Theorem~\ref{th2}, $f_2(x)\in {\cA\cH}_0(\bR)$.

3. Let $f_k(x)$ be a QAHD of degree $\lambda$
and of order $k$. Let us assume that any QAHD
$f_j(x)$, \ $j=0,1,\dots,k-1$ is a distribution of degree $\lambda$
and of order $j$ (in the sense of Definition~\ref{de5}).
Then, according to Theorem~\ref{th2}, $f_j(x)\in {\cA\cH}_0(\bR)$ and
relation (\ref{17.1}) holds.
Thus, in view of our assumption, (\ref{50}) and (\ref{17.1}) imply
for all $a,b>0$:
$$
f_k(abx)=(ab)^{\lambda}f_k(x)+\sum_{r=1}^{k}h_r(ab)f_{k-r}(x)
=a^{\lambda}f_k(bx)+\sum_{r=1}^{k}h_r(a)f_{k-r}(bx)
$$
$$
=a^{\lambda}\Big(b^{\lambda}f_k(x)+\sum_{r=1}^{k}h_r(b)f_{k-r}(x)\Big)
\qquad\qquad\qquad\qquad\qquad
$$
$$
\qquad\qquad\qquad
+\sum_{r=1}^{k-1}h_r(a)\Big(b^{\lambda}f_{k-r}(x)
+\sum_{j=1}^{k-r}b^{\lambda}\log^j{b}f^{(k-r)}_{k-r-j}(x)\Big)
+h_k(a)b^{\lambda}f_0(x)
$$
$$
=(ab)^{\lambda}f_k(x)
+\sum_{r=1}^{k}\Big(a^{\lambda}h_{r}(b)+b^{\lambda}h_{r}(a)\Big)f_{k-r}(x)
\qquad\qquad
$$
$$
\qquad\qquad\qquad\qquad\qquad
+\sum_{r=1}^{k-1}\sum_{j=1}^{k-r}
h_r(a)b^{\lambda}\log^j{b} f^{(k-r)}_{k-r-j}(x),
$$
where $f^{(k-r)}_{k-r-j}(x)$ is a distribution of degree $\lambda$
and of order $k-r-j$ (in the sense of Definition~\ref{de5}) which
belongs to ${\cA\cH}_0(\bR)$, \ $r=1,\dots,k-1$, \
$j=1,\dots,k-r$. By changing the sum order, one can easily see that
for all $a,b>0$:
$$
\sum_{r=1}^{k}h_r(ab)f_{k-r}(x)
=\sum_{r=1}^{k}\Big(a^{\lambda}h_{r}(b)+b^{\lambda}h_{r}(a)\Big)f_{k-r}(x)
\qquad\qquad\qquad\qquad
$$
\begin{equation}
\label{59}
\qquad\qquad
+\sum_{r=2}^{k}\sum_{j=1}^{r-1}
h_{r-j}(a)b^{\lambda}\log^j{b} f^{(k-r+j)}_{k-r}(x).
\end{equation}

Since, in view of Lemma~\ref{lem2}, a distribution
$f_{k-1}\in {\cA\cH}_1(\bR)_0$ of order $k-1$ and distributions
$f_{k-r},f^{(k-r+j)}_{k-r}\in {\cA\cH}_1(\bR)_0$ of order $k-r$,
are linear independent, $r=2,\dots,k$, \ $j=1,\dots,r-1$,
relation (\ref{59}) implies that for all $a,b>0$
$$
\begin{array}{rcl}
\displaystyle
h_1(ab)&=&a^{\lambda}h_1(b)+b^{\lambda}h_1(a), \medskip \\
\displaystyle
h_2(ab)f_{k-2}&=&\Big(a^{\lambda}h_2(b)+b^{\lambda}h_2(a)\Big)f_{k-2}
+h_1(a)b^{\lambda}\log b f^{(k-1)}_{k-2}, \medskip \\
\displaystyle
h_3(ab)f_{k-3}&=&\Big(a^{\lambda}h_3(b)+b^{\lambda}h_3(a)\Big)f_{k-3}
+\sum_{j=1}^{2}h_{3-j}(a)b^{\lambda}\log^j b f^{(k-3+j)}_{k-3}, \medskip  \\
\displaystyle
\cdots\cdots\cdots&\cdot&\cdots\cdots\cdots\cdots\cdots\cdots\cdots
\cdots\cdots\cdots\cdots\cdots\cdots\cdots\cdots\cdots\cdots, \medskip  \\
\displaystyle
h_k(ab)f_{0}&=&\Big(a^{\lambda}h_k(b)+b^{\lambda}h_k(a)\Big)f_{0}
+\sum_{j=1}^{k-1}h_{k-j}(a)b^{\lambda}\log^j b f^{(j)}_{0}. \\
\end{array}
$$

Taking into account that the function $\frac{h_j(ab)}{(ab)^{\lambda}}
-\frac{h_j(a)}{(a)^{\lambda}}-\frac{h_j(b)}{(b)^{\lambda}}$ is symmetric
in $a$ and $b$, it is easy to see that the last system has a {\it non-trivial}
solution only if $f^{(k-r+j)}_{k-r}(x)=C^{(k-r+j)}_{k-r}f_{k-r}(x)$,
where $C^{(k-r+j)}_{k-r}$ are constants, $r=2,3,\dots,k$, \
$j=1,2,\dots,r-1$. Thus in view of Lemma~\ref{lem2}, we obtain the following
system of functional equations
\begin{equation}
\label{60}
\begin{array}{rcl}
\displaystyle
h_1(ab)&=&a^{\lambda}h_1(b)+b^{\lambda}h_1(a), \medskip  \\
\displaystyle
h_2(ab)&=&a^{\lambda}h_2(b)+b^{\lambda}h_2(a)
+C^{(k-1)}_{k-2}h_1(a)b^{\lambda}\log b, \medskip  \\
\displaystyle
h_3(ab)&=&a^{\lambda}h_3(b)+b^{\lambda}h_3(a)
+\sum_{j=1}^{2}C^{(k-3+j)}_{k-3}h_{3-j}(a)b^{\lambda}\log^j b, \medskip  \\
\displaystyle
\cdots\cdots&\cdot&\cdots\cdots\cdot\cdots\cdots\cdots
\cdots\cdots\cdots\cdots\cdots\cdots\cdots\cdots\cdots\cdots, \medskip  \\
\displaystyle
h_k(ab)&=&a^{\lambda}h_k(b)+b^{\lambda}h_k(a)
+\sum_{j=1}^{k-1}C^{(j)}_{0}h_{k-j}(a)b^{\lambda}\log^j b. \\
\end{array}
\end{equation}
Consequently, the functions
${\widetilde h}_j(a)=\frac{h_j(a)}{a^{\lambda}}$ satisfy the system of equation
\begin{equation}
\label{60.1}
\begin{array}{rcl}
\displaystyle
{\widetilde h}_1(ab)&=&{\widetilde h}_1(b)+{\widetilde h}_1(a),
\medskip  \\
\displaystyle
{\widetilde h}_2(ab)&=&{\widetilde h}_2(b)+{\widetilde h}_2(a)
+C^{(k-1)}_{k-2}{\widetilde h}_1(a)\log b, \medskip  \\
\displaystyle
{\widetilde h}_3(ab)&=&{\widetilde h}_3(b)+{\widetilde h}_3(a)
+\sum_{j=1}^{2}C^{(k-3+j)}_{k-3}{\widetilde h}_{3-j}(a)\log^j b, \medskip  \\
\displaystyle
\cdots\cdots&\cdot&\cdots\cdots\cdot\cdots\cdots\cdots
\cdots\cdots\cdots\cdots\cdots\cdots\cdots\cdots, \medskip  \\
\displaystyle
{\widetilde h}_k(ab)&=&{\widetilde h}_k(b)+{\widetilde h}_k(a)
+\sum_{j=1}^{k-1}C^{(j)}_{0}{\widetilde h}_{k-j}(a)\log^j b, \\
\end{array}
\end{equation}
where ${\widetilde h}_j(1)=0$, $j=1,2,\dots,k$.

By changing variables $\psi_j(z)={\widetilde h}_j(e^{z})$, where
$\psi_j(0)=0$, $j=1,2,\dots,k$ and $a=e^{\xi}$, $b=e^{\eta}$, system
(\ref{60.1}) can be rewritten as
\begin{equation}
\label{60.2}
\begin{array}{rcl}
\displaystyle
\psi_1(\xi+\eta)&=&\psi_1(\xi)+\psi_1(\eta),
\medskip  \\
\displaystyle
\psi_2(\xi+\eta)&=&\psi_2(\xi)+\psi_2(\eta)
+C^{(k-1)}_{k-2}\psi_1(\xi)\eta, \medskip  \\
\displaystyle
\psi_3(\xi+\eta)&=&\psi_3(\xi)+\psi_3(\eta)
+\sum_{j=1}^{2}C^{(k-3+j)}_{k-3}\psi_{3-j}(\xi)\eta^j, \medskip  \\
\displaystyle
\cdots\cdots\cdots&\cdot&\cdots\cdots\cdot\cdots\cdots
\cdots\cdots\cdots\cdots\cdots\cdots\cdots\cdots, \medskip  \\
\displaystyle
\psi_k(\xi+\eta)&=&\psi_k(\xi)+\psi_k(\eta)
+\sum_{j=1}^{k-1}C^{(j)}_{0}\psi_{k-j}(\xi)\eta^j. \\
\end{array}
\end{equation}
Differentiating relations (\ref{60.2}) with respect to $\eta$
and setting  $\eta=0$, we obtain a system of differential equations
\begin{equation}
\label{60.3}
\begin{array}{rcl}
\displaystyle
\psi_1'(\xi)&=&\psi_1'(0),
\medskip  \\
\displaystyle
\psi_2'(\xi)&=&\psi_2'(0)
+C^{(k-1)}_{k-2}\psi_1(\xi), \medskip  \\
\displaystyle
\psi_3'(\xi)&=&\psi_3'(0)+C^{(k-3+j)}_{k-1}\psi_{2}(\xi), \medskip  \\
\displaystyle
\cdots\cdots&\cdot&\cdots\cdots\cdot\cdots\cdots\cdots\cdots, \medskip  \\
\displaystyle
\psi_k'(\xi)&=&\psi_k'(0)+C^{(1)}_{0}\psi_{k-1}(\xi), \\
\end{array}
\end{equation}
where $\psi_j(0)=0$, $j=1,2,\dots,k$.

By successive integration it is easy to see that a solution of
system (\ref{60.3}) has the form
$$
\psi_r(\xi)=\sum_{j=1}^{r}A_{r}^{j}\xi^j,
$$
where $A_r^j$ are constants, which can be calculated, $r=1,2,\dots,k$,
$j=1,2,\dots,r$.

Since $a=e^{\xi}$, \ $\psi_j(z)={\widetilde h}_j(e^{z})$, then
${\widetilde h}_r(a)=\sum_{j=1}^{r}A_{r}^{j}\log^j a$ and
\begin{equation}
\label{61}
h_r(a)=a^{\lambda}\sum_{j=1}^{r}A_r^j \log^j a,
\end{equation}
where $A_r^j$ are constant, $r=1,2,\dots,k$, $j=1,2,\dots,r$.

By substituting functions (\ref{61}) into relation (\ref{50}), the last
relation can be rewritten in the form (\ref{50}), i.e., as
$$
f_k(ax)=a^{\lambda}f_k(x)
+a^{\lambda}\sum_{r=1}^{k}\sum_{j=1}^{r}A_r^j \log^j a f_{k-r}(x)
$$
\begin{equation}
\label{62}
=a^{\lambda}f_k(x)
+\sum_{r=1}^{k}a^{\lambda}\log^r{a}{\widetilde f}_{k-r}(x),
\end{equation}
where according to our assumption, distribution
${\widetilde f}_{k-r}(x)=\sum_{j=r}^{k}A_{j}^{r} f_{k-j}(x)$
belongs to the class ${\cA\cH}_0(\bR)$, \ $r=1,2,\dots,k$. Moreover, in
view of Remark~\ref{rem2}, ${\widetilde f}_{k-r}(x)$ is a distribution
of degree $\lambda$ and of order $k-r$ (in the sense of Definition~\ref{de5}),
$r=1,2,\dots,k$. Consequently, $f_k(x)$ satisfies relation (\ref{17.1}).

Thus, according of the induction axiom, the theorem is proved.
\end{proof}

\subsection{Resume.}\label{s4.3}
In~\cite[Ch.I,\S 4.1.]{G-Sh} it was proved that in order to introduce an AHD
of order $k=1$, one can use Definition~\ref{de3} instead of Definition~(\ref{1.1}).
Similarly, according to Theorem~\ref{th3}, in order to introduce a QAHD,
instead of Definition~\ref{de6} one can use the following definition
(in fact, Definition~\ref{de5}).

\begin{definition}
\label{de7} \rm
A distribution $f_{k}\in {\cD}'(\bR)$ is called a {\it QAHD} of
degree~$\lambda$ and of order $k$, $k=0,1,2,\dots$,
if for any $a >0$ and $\varphi \in {\cD}(\bR)$
\begin{equation}
\label{17**}
\Bigl\langle f_{k}(x),\varphi\Big(\frac{x}{a}\Big) \Bigr\rangle
=a^{\lambda+1} \bigl\langle f_{k}(x),\varphi(x) \bigr\rangle
+\sum_{r=1}^{k}
a^{\lambda+1}\log^r{a}\bigl\langle f_{k-r}(x),\varphi(x) \bigr\rangle,
\end{equation}
i.e.,
\begin{equation}
\label{17.1**}
f_{k}(ax)=a^{\lambda}f_{k}(x)+\sum_{r=1}^{k}a^{\lambda}\log^r{a} f_{k-r}(x),
\end{equation}
where $f_{k-r}(x)$ is an AHD of degree $\lambda$ and of order $k-r$,
\ $r=1,2,\dots,k$. Here for $k=0$ we suppose that the
sums in the right-hand sides of (\ref{17**}), (\ref{17.1**}) are empty.
\end{definition}

Thus instead of the term {\it ``distribution of degree~$\lambda$
and of order $k$''} one can use the term {\it ``QAHD of degree~$\lambda$
and of order $k$''}.

According to Remark~\ref{rem2}, the sum of a QAHD
of degree~$\lambda$ and of order~$k$, and a QAHD
of degree~$\lambda$ and of order~$r\le k-1$
is a QAHD of degree~$\lambda$ and of order~$k$.

According to Theorems~\ref{th2},~\ref{th3}, the class of QAHDs
coincides with the Gel$'$fand--Shilov class ${\cA\cH}_0(\bR)$.

\section{Multidimensional QAHDs}
\label{s5}

\begin{definition}
\label{de8} \rm
(see~\cite[Ch.III,\S 3.1.,(1)]{G-Sh})
A distribution $f_{0}(x)=f_{0}(x_1,\dots,x_n)$ from ${\cD}'(\bR^n)$ is
called {\it homogeneous} of degree~$\lambda$ if for any $a >0$ and
$\varphi \in {\cD}(\bR^n)$
$$
\Bigl\langle f_{0},\varphi\Big(\frac{x_1}{a},\dots,\frac{x_n}{a}\Big)\Bigr\rangle
=a^{\lambda+n} \bigl\langle f_{0},\varphi(x_1,\dots,x_n) \bigr\rangle
$$
i.e.,
$$
f_{0}(ax_1,\dots,ax_n)=a^{\lambda}f_{0}(x_1,\dots,x_n).
$$
\end{definition}

Recall a well-known theorem.
\begin{theorem}
\label{th5}
{\rm(}see~{\rm\cite[Ch.III,\S 3.1.]{G-Sh}}{\rm)}
A distribution $f_{0}(x)$ is homogeneous of degree~$\lambda$ if and only if
it satisfies the Euler equation
$$
\sum_{j=1}^{n}x_j\frac{\partial f_{0}}{\partial x_j}=\lambda f_{0}.
$$
\end{theorem}

Now we introduce a multidimensional analog of Definition~\ref{de7} and
prove a multidimensional analog of Theorem~\ref{th5}.
\begin{definition}
\label{de9} \rm
We say that a distribution $f_{k} \in {\cD}'(\bR^n)$ is a {\it QAHD\/}
of degree $\lambda$ and of order~$k$, \ $k=0,1,2,\dots$, if for any
$a >0$ we have
\begin{equation}
\label{17.1***}
f_{k}(ax)=f_k(ax_1,\dots,ax_n)
=a^{\lambda}f_{k}(x)+\sum_{r=1}^{k}a^{\lambda}\log^r{a} f_{k-r}(x),
\end{equation}
where $f_{k-r}(x)$ is a QAHD of degree $\lambda$ and
of order $k-r$, \ $r=1,2,\dots,k$. Here for $k=0$ we suppose that the
sum in the right-hand side of (\ref{17.1}) is empty.
\end{definition}

\begin{theorem}
\label{th6}
$f_k(x)$ is a QAHD of degree~$\lambda$ and of order $k$, $k\ge 1$
if and only if it satisfies the Euler type system of equations,
i.e, there exist distributions $f_{k-1},\dots,f_{0}$ such that
\begin{equation}
\label{75}
\begin{array}{rcl}
\displaystyle
\sum_{j=1}^{n}x_j\frac{\partial f_{k}}{\partial x_j}&=&\lambda f_{k}+f_{k-1}, \\
\displaystyle
\sum_{j=1}^{n}x_j\frac{\partial f_{k-1}}{\partial x_j}&=&\lambda f_{k-1}+f_{k-2}, \\
\displaystyle
\cdots\cdots\cdots &\cdot&\cdots\cdots\cdots, \\
\displaystyle
\sum_{j=1}^{n}x_j\frac{\partial f_{1}}{\partial x_j}&=&\lambda f_{1}+f_{0}, \\
\displaystyle
\sum_{j=1}^{n}x_j\frac{\partial f_{0}}{\partial x_j}&=&\lambda f_{0}, \\
\end{array}
\end{equation}
i.e., for all $\varphi \in {\cD}(\bR^n)$
$$
\begin{array}{rcl}
\displaystyle
-\Bigl\langle f_{k},
\sum_{j=1}^{n}x_j\frac{\partial \varphi}{\partial x_j}\Bigr\rangle
&=&(\lambda+n)\langle f_{k},\varphi\rangle+\langle f_{k-1},\varphi\rangle, \\
\displaystyle
-\Bigl\langle f_{k-1},
\sum_{j=1}^{n}x_j\frac{\partial \varphi}{\partial x_j}\Bigr\rangle
&=&(\lambda+n)\langle f_{k-1},\varphi\rangle+\langle f_{k-2},\varphi\rangle, \\
\displaystyle
\cdots\cdots\cdots\cdots\cdots &\cdot&\cdots\cdots\cdots\cdots\cdots
\cdots\cdots\cdots, \\
\displaystyle
-\Bigl\langle f_{1},
\sum_{j=1}^{n}x_j\frac{\partial \varphi}{\partial x_j}\Bigr\rangle
&=&(\lambda+n)\langle f_{1},\varphi\rangle+\langle f_{0},\varphi\rangle, \\
\displaystyle
-\Bigl\langle f_{0},
\sum_{j=1}^{n}x_j\frac{\partial \varphi}{\partial x_j}\Bigr\rangle
&=&(\lambda+n)\langle f_{0},\varphi\rangle, \\
\end{array}
$$
\end{theorem}

\begin{proof}
Let $f_{k}\in {\cD}'(\bR^n)$ be a QAHD of degree $\lambda$ and of
order~$k$. This implies that, according to Definition~\ref{de9},
there are distributions $f_{j}$, $j=0,1,2,\dots,k-1$ and
$f_{k-s-r}^{(k-s)}$, $s=0,1,2,\dots,k-2$, $r=2,\dots,k-s$ such that
\begin{equation}
\label{76}
\begin{array}{rcl}
\displaystyle
f_k(ax_1,\dots,ax_n)&=&a^{\lambda}f_{k}(x)+a^{\lambda}\log{a} f_{k-1}(x)
+\sum_{r=2}^{k}a^{\lambda}\log^r{a} f_{k-r}^{(k)}(x), \medskip \\
\displaystyle
f_{k-1}(ax_1,\dots,ax_n)&=&a^{\lambda}f_{k-1}(x)+a^{\lambda}\log{a} f_{k-2}(x)
+\sum_{r=2}^{k-1}a^{\lambda}\log^r{a} f_{k-1-r}^{(k-1)}(x), \medskip \\
\displaystyle
f_{k-2}(ax_1,\dots,ax_n)&=&a^{\lambda}f_{k-2}(x)+a^{\lambda}\log{a} f_{k-3}(x)
+\sum_{r=2}^{k-2}a^{\lambda}\log^r{a} f_{k-2-r}^{(k-2)}(x), \medskip \\
\displaystyle
\cdots\cdots\cdots\cdots\cdots &\cdot&\cdots\cdots\cdots\cdots\cdots\cdots
\cdots\cdots\cdots\cdots\cdots\cdots\cdots\cdots\cdots\cdots, \medskip \\
\displaystyle
f_1(ax_1,\dots,ax_n)&=&a^{\lambda}f_{1}(x)+a^{\lambda}\log{a} f_{0}(x), \medskip \\
\displaystyle
f_0(ax_1,\dots,ax_n)&=&a^{\lambda}f_{0}(x). \\
\end{array}
\end{equation}
Differentiating (\ref{76}) with respect to $a$ and setting $a=1$, we obtain
system (\ref{75}).

Conversely, let $f_{k}\in {\cD}'(\bR^n)$ be a distribution satisfying
system (\ref{75}), i.e., there are distributions $f_{j}\in {\cD}'(\bR^n)$,
$j=0,1,2,\dots,k-1$ such that system (\ref{75}) holds. We prove by
induction that $f_{k}$ is a QAHD of degree $\lambda$ and of order~$k$.

Let $k=0$. This fact is proved by Theorem~\ref{th5}.

If $k=1$ then the following system of equations
\begin{equation}
\label{76.1}
\begin{array}{rcl}
\displaystyle
\sum_{j=1}^{n}x_j\frac{\partial f_{1}}{\partial x_j}&=&\lambda f_{1}+f_{0}, \\
\displaystyle
\sum_{j=1}^{n}x_j\frac{\partial f_{0}}{\partial x_j}&=&\lambda f_{0} \\
\end{array}
\end{equation}
holds. Here, in view of Theorem~\ref{th5} the second equation in
(\ref{76.1}) implies that $f_{0}$ is a HD.

Consider the function
$$
g_1(a)=f_1(ax_1,\dots,ax_n)-a^{\lambda}f_{1}(x)-a^{\lambda}\log{a} f_{0}(x).
$$
It is clear that $g_1(1)=0$. By differentiation
we have
\begin{equation}
\label{76.2}
g_1'(a)=\sum_{j=1}^{n}x_j\frac{\partial f_1}{\partial x_j}(ax_1,\dots,ax_n)
-\lambda a^{\lambda-1}f_{1}(x)-(\lambda a^{\lambda-1}\log{a}+a^{\lambda-1})f_{0}(x)
\end{equation}
Applying the first relation in (\ref{76.1}) to the arguments $ax_1,\dots,ax_n$
we find that
\begin{equation}
\label{76.3}
\sum_{j=1}^{n}x_j\frac{\partial f_{1}}{\partial x_j}(ax_1,\dots,ax_n)=
\frac{\lambda}{a}f_{1}(ax_1,\dots,ax_n)+\frac{1}{a}f_{0}(ax_1,\dots,ax_n). \\
\end{equation}
Substituting (\ref{76.3}) into (\ref{76.2}) and taking into account that
$\frac{1}{a}f_{0}(ax_1,\dots,ax_n)=a^{\lambda-1}f_{0}$, we find that $g_1(a)$
satisfies the differential equation
\begin{equation}
\label{76.4}
g_1'(a)=\frac{\lambda}{a}g_1(a), \qquad g_1(1)=0.
\end{equation}
Obviously, its solution is $g_1(a)g_1(a)=0$. Thus
$g_1(a)=f_1(ax_1,\dots,ax_n)-a^{\lambda}f_{1}(x)
-a^{\lambda}\log{a} f_{0}(x)=0$, i.e, $f_1(x)$ is an AHD of order $k=1$,
i.e., a QAHD of order $k=1$.

Let us assume that for $k-1$ the theorem holds, i.e., if $f_{k-1}$
satisfies all the equations in (\ref{75}) except the first one, then
$f_{k-1}$ is a QAHD of degree $\lambda$ and of order~$k-1$.

Let us consider the case $k$. Let the Euler type system (\ref{75})
be satisfied, i.e, there exist distributions $f_{k-1}, \dots, f_{0}$
such that (\ref{75}) holds. Note that in view of our assumption,
$f_{k-1}$ is a QAHD of order $k-1$.

Consider the function
\begin{equation}
\label{76.5*}
g_k(a)=f_k(ax_1,\dots,ax_n)-a^{\lambda}f_{k}(x)-a^{\lambda}\log{a} f_{k-1}(x).
\end{equation}
It is clear that $g_k(1)=0$. By differentiation we have
\begin{equation}
\label{76.5}
g_k'(a)=\sum_{j=1}^{n}x_j\frac{\partial f_k}{\partial x_j}(ax_1,\dots,ax_n)
-\lambda a^{\lambda-1}f_{k}(x)-(\lambda a^{\lambda-1}\log{a}+a^{\lambda-1})f_{k-1}(x)
\end{equation}
Applying the first relation in (\ref{75}) to the arguments $ax_1,\dots,ax_n$
we find that
\begin{equation}
\label{76.6}
\sum_{j=1}^{n}x_j\frac{\partial f_{k}}{\partial x_j}(ax_1,\dots,ax_n)=
\frac{\lambda}{a}f_{k}(ax_1,\dots,ax_n)+\frac{1}{a}f_{k-1}(ax_1,\dots,ax_n). \\
\end{equation}
Substituting (\ref{76.6}) into (\ref{76.5}) and taking into account that
(in view of our assumption) $f_{k-1}$ is a QAHD of order $k-1$, i.e.,
$$
f_{k-1}(ax_1,\dots,ax_n)=a^{\lambda}f_{k-1}(x)
+\sum_{r=1}^{k-1}a^{\lambda}\log^r{a} f_{k-1-r}^{(k-1)}(x),
$$
where $f_{k-1-r}^{(k-1)}(x)$ is a QAHD of order $k-1-r$, \ $r=1,2,\dots,k-1$,
we find that $g_k(a)$ satisfies the linear differential equation
\begin{equation}
\label{76.7}
g_k'(a)=\frac{\lambda}{a}g_k(a)
+\sum_{r=1}^{k-1}a^{\lambda-1}\log^r{a} f_{k-1-r}^{(k-1)}(x), \qquad g_1(1)=0.
\end{equation}
Now it is easy to see that its general solution has the form
$$
g_k(a)=\sum_{r=1}^{k-1}a^{\lambda}\log^{r+1}{a}\frac{ f_{k-1-r}^{(k-1)}(x)}{r+1}
+a^{\lambda}C(x),
$$
where $C(x)$ is a distribution. Taking into account that $g_1(1)=0$, we
calculate $C(x)=0$. Thus
\begin{equation}
\label{76.8}
g_k(a)=\sum_{r=1}^{k-1}a^{\lambda}\log^{r+1}{a}\frac{ f_{k-1-r}^{(k-1)}(x)}{r+1}.
\end{equation}

By substituting (\ref{76.8}) into (\ref{76.5*}), we find
\begin{equation}
\label{76.9}
f_k(ax_1,\dots,ax_n)=a^{\lambda}f_{k}(x)-a^{\lambda}\log{a} f_{k-1}(x)
+\sum_{r=2}^{k}a^{\lambda}\log^{r}{a}\frac{ f_{k-r}^{(k-1)}(x)}{r},
\end{equation}
where by our assumption $f_{k-1}$ is a QAHD of order $k-1$, and,
consequently, $f_{k-r}^{(k-1)}(x)$ is a QAHD of order $k-r$, \ $r=2,\dots,k$.
Thus, in view of Definition~\ref{de9}, $f_k$ is a QAHD of order $k$.

Thus, according to the induction axiom, the theorem is proved.
\end{proof}

\section{The Fourier transform of QAHDs}
\label{s6}

\subsection{The Fourier transform.}
\label{s6.1}
The Fourier transform of $\varphi\in {\cD}(\bR^n)$ is defined as
$$
F[\varphi](\xi)=\int_{\bR^n}\varphi(x)e^{i\xi\cdot x}\,d^nx,
\quad \xi \in \bR^n,
$$
where $\xi\cdot x$ is the scalar product of vectors $x$ and $\xi$.
We define the Fourier transform $F[f]$ of a distribution~\cite[Ch.II]{G-Sh}
$$
\langle F[f],\varphi\rangle=\langle f,F[\varphi]\rangle,
\quad \forall \, \varphi\in {\cD}(\bR^n).
$$
Let $f\in{\cD}'(\bR^n)$. If $a\ne 0$ is a constant then
\begin{equation}
\label{76*}
F[f(ax)](\xi)=F[f(ax_1,\dots,ax_n)](\xi)=|a|^{-n}F[f(x)]\Big(\frac{\xi}{a}\Big).
\end{equation}

\begin{theorem}
\label{th7}
If $f\in {\cD}'(\bR^n)$ is a QAHD of degree $\lambda$ and of order $k$,
then its Fourier transform $F[f]$ is a QAHD of degree $-\lambda-1$
and of order $k$, $k=0,1,2,\dots$.
\end{theorem}

\begin{proof}
We prove this theorem by induction.

If $k=0$ then using (\ref{76*}) and Definition~\ref{de8},
we have for all $a>0$
\begin{equation}
\label{77}
F\big[f(x)\big](a\xi)=a^{-n}F\big[f\Big(\frac{x}{a}\Big)\big](\xi)
=a^{-\lambda-n}F\big[f(x)\big](\xi),
\end{equation}
i.e., $F[f(x)](\xi)$ is a HD of degree $-\lambda-n$.

Let $k=1$. Using (\ref{76*}) and Definition~\ref{de9},
we obtain for all $a>0$
$$
F\big[f(x)\big](a\xi)=a^{-n}F\big[f\Big(\frac{x}{a}\Big)\big](\xi)
\qquad\qquad\qquad\qquad\qquad\qquad\quad
$$
$$
\qquad\qquad
=a^{-\lambda-n}F\big[f(x)\big](\xi)-a^{-\lambda-n}\log{a}F\big[f_0(x)\big](\xi),
$$
where $f_{0}$ is a HD of degree $\lambda$.
In view of (\ref{77}), $F[f_0](\xi)$ is a HD of degree $-\lambda-n$,
hence, according to Definition~\ref{de9}, $F[f(x)](\xi)$ is an AHD of
degree $-\lambda-n$ and of order $k=1$, i.e., a QAHD of
degree $-\lambda-n$ and of order $k=1$.

Let $f$ be a QAHD of degree $\lambda$ and order $k$, \ $k=2,3,\dots$.
By using (\ref{76*}) and Definition~\ref{de9}, for all $a>0$ we have
$$
F\big[f(x)\big](a\xi)=a^{-n}F\big[f\Big(\frac{x}{a}\Big)\big](\xi)
\qquad\qquad\qquad\qquad\qquad\qquad\qquad\qquad\qquad\qquad
$$
$$
\qquad
=a^{-\lambda-n}F\big[f(x)\big](\xi)
+\sum_{r=1}^{k}(-1)^ra^{-\lambda-n}\log^r{a}F\big[f_{k-r}(x)\big](\xi),
$$
where $f_{k-r}(x)$ is a QAHD of degree~$\lambda$ and order $k-r$,
\ $r=1,2,\dots,k$.

Suppose that the theorem holds for QAHDs of degree~$\lambda$ and order
$k=1,2\dots,k-1$. Hence, by induction the last relation implies that
$F[f](\xi)$ is a QAHD of degree~$\lambda$ and of order $k$.

The theorem is thus proved.

Taking into account Theorem~\ref{th2} and Remark~\ref{rem2},
one can prove this theorem directly by calculating the Fourier transform
of distributions $P\big(x_{\pm}^{-n}\log^{k-1} x_{\pm}\big)$
and $x_{\pm}^{\lambda}\log^k x_{\pm}$, where $\lambda \ne -1,-2,\dots$.
\end{proof}

Thus $F\big[{\cA\cH}_0(\bR)\big]={\cA\cH}_0(\bR)$.

\subsection{Gamma functions generated by QAHDs.}\label{s6.2}
Consider the Fourier transform of homogeneous distribution
$x_{+}^{\lambda}$, $\lambda \ne -1,-2,\dots$, which according
to Theorem~\ref{th7}, is represented as
\begin{equation}
\label{80}
F\big[x_{+}^{\lambda}](\xi)=C(\xi+i0)^{-\lambda-1},
\end{equation}
where $C$ is a constant and the distribution $(x\pm i0)^{\lambda}$
is given by (\ref{8.1}). Setting $\xi=i$, one can calculate that
\begin{equation}
\label{80.1}
C=i^{\lambda+1}\int_{0}^{\infty}x^{\lambda}e^{-x}\,dx
=i^{\lambda+1}\Gamma(\lambda+1).
\end{equation}
Thus the factor of proportionality in (\ref{80}) is (up to $i^{\lambda+1}$)
the $\Gamma$-function, $\Gamma(\lambda+1)=\int_{0}^{\infty}x^{\lambda}e^{-x}\,dx$.

In view of Theorem~\ref{th7} and Remark~\ref{rem2},
we have for $\lambda \ne -1,-2,\dots$
\begin{equation}
\label{81}
F\big[x_{+}^{\lambda}\log^k{x_{+}}\big](\xi)
=\sum_{j=0}^{k}A_{k-j}(\xi+i0)^{-\lambda-1}\log^{k-j}(\xi+i0),
\end{equation}
and for  $\lambda=-n$, $n\in \bN$
\begin{equation}
\label{82}
F\big[P\big(x_{+}^{-n}\log^{k-1}{x_{+}}\big)\big](\xi)
=\sum_{j=0}^{k}B_{k-j}\xi^{n-1}\log^{k-j}(\xi+i0),
\end{equation}
where $A_j$, $B_j$ are constants, $j=1,\dots,k$.
Here $(\xi+i0)^{-\lambda-1}\log^{k-j}(\xi+i0)$ and
$\xi^{n-1}\log^{k-j}(\xi+i0$ are QAHDs of order $k-j$ and
of degree~$-\lambda-1$ and $n-1$, respectively (we set
$(\xi+i0)^{n-1}\equiv \xi^{n-1}$, $n\in \bN$).

Similarly (\ref{80.1}), we call the factors
\begin{equation}
\label{82.1}
\Gamma_j(\lambda+1;k)=i^{-\lambda-1}\log^{j}{i}\,A_j,
\end{equation}
and
\begin{equation}
\label{82.2}
\Gamma_j(-n+1;k)=i^{n-1}\log^{j}{i}\,B_j
\end{equation}
the {\it associated homogeneous $j-\Gamma$-functions} of order $k$ and of
degree $\lambda$ ($\lambda \ne -n$) and $-n$, respectively,
$j=0,1,\dots,k$,\ $n\in \bN$.

By successive substituting $\xi=i,2i,\dots,(k+1)i$ into
(\ref{81}) and (\ref{82}), we obtain a linear system of equation for
$A_0,\dots,A_k$ and $B_0,\dots,B_k$. Solving these systems, one can
calculate {\it associated homogeneous $\Gamma$-functions}
$\Gamma_j(\lambda+1;k)$ and $\Gamma_j(-n+1;k)$, respectively.

Now we calculate the $\Gamma$-functions in particular case $k=1$.

Let $\lambda \ne -1,-2,\dots$. According to~\cite[Ch.II,\S2.4.,(1)]{G-Sh},
$$
F\big[x_{+}^{\lambda}\log{x_{+}}\big](\xi)
=-i^{\lambda+1}\Gamma(\lambda+1)(\xi+i0)^{-\lambda-1}\log(\xi+i0)
\qquad\qquad\qquad
$$
$$
\qquad\qquad\quad
+i^{\lambda+1}\Big(\Gamma'(\lambda+1)+i\frac{\pi}{2}\Gamma(\lambda+1)\Big)
(\xi+i0)^{-\lambda-1}.
$$
This relation and (\ref{82.1}) imply that
\begin{equation}
\label{83}
\begin{array}{rcl}
\displaystyle
\Gamma_1(\lambda+1;1)&=&-i\frac{\pi}{2}\,\Gamma(\lambda+1), \medskip \\
\displaystyle
\Gamma_0(\lambda+1;1)&=&\Gamma'(\lambda+1)+i\frac{\pi}{2}\Gamma(\lambda+1).
\end{array}
\end{equation}

Let $\lambda=-1,-2,\dots$. According to~\cite[Ch.II,\S2.4.,(14)]{G-Sh},
$$
F\big[x_{+}^{-n}\big](\xi)
=-a^{(n)}_{-1}\xi^{n-1}\log(\xi+i0)+a^{(n)}_{0}\xi^{n-1},
$$
where
$$
\begin{array}{rcl}
\displaystyle
a^{(n)}_{-1}&=&\frac{i^{n+1}}{(n-1)!}, \medskip \\
\displaystyle
a^{(n)}_{0}&=&\frac{i^{n+1}}{(n-1)!}
\Big(1+\frac{1}{2}+\cdots+\frac{1}{n-1}+\Gamma'(1)+i\frac{\pi}{2}\Big).
\end{array}
$$
Thus, in view of (\ref{82.2}), we have
\begin{equation}
\label{85}
\begin{array}{rcl}
\displaystyle
\Gamma_1(-n+1;1)&=&-i\frac{\pi}{2}\frac{(-1)^n}{(n-1)!}, \medskip \\
\displaystyle
\Gamma_0(-n+1;1)&=&\frac{(-1)^n}{(n-1)!}
\Big(1+\frac{1}{2}+\cdots+\frac{1}{n-1}+\Gamma'(1)+i\frac{\pi}{2}\Big).
\end{array}
\end{equation}

Of course, formulas (\ref{83}), (\ref{85}) can be derived directly.

According to (\ref{83}), (\ref{85}), we have
$$
\begin{array}{rcl}
\displaystyle
\Gamma_1(\lambda+1;1)&=&\lambda\Gamma_1(\lambda;1), \medskip \\
\displaystyle
\Gamma_0(\lambda+1;1)&=&\lambda\Gamma_0(\lambda;1)+\Gamma(\lambda);
\end{array}
$$
and
$$
\begin{array}{rcl}
\displaystyle
\Gamma_1(-n+1;1)&=&(-n)\Gamma_1(-n;1), \\
\displaystyle
\Gamma_0(-n+1;1)&=&(-n)\Gamma_0(-n;1)-\frac{(-1)^n}{n!},
\end{array}
$$
where ${\rm res}_{\lambda=-n}\Gamma(\lambda)=\frac{(-1)^n}{n!}$.

The other {\it associated homogeneous $\Gamma$-functions}
$\Gamma_j(\lambda+1;k)$ and $\Gamma_j(-n+1;k)$ can be calculated
in the same way and their properties can be studied. But here we
omit these problems.

\begin{center}
{\bf Acknowledgements}
\end{center}

The author is greatly indebted to V.~I.~Polischook for many fruitful
discussions and careful reading all versions of this paper.
The author is also grateful to Evgenia Sorokina whose remarks helped
to prove Theorem~\ref{th3}.

\end{document}